\documentclass[preprint,12pt]{elsarticle}
\biboptions{sort&compress}

\usepackage{etoolbox}
\patchcmd{\keyword}{\textit{Keywords:}}{\textbf{Keywords:}}{}{}

\usepackage{amssymb}
\usepackage{amsmath}
\usepackage{booktabs}
\usepackage{float}
\usepackage{subcaption}
\usepackage{hyperref}

\journal{Nuclear Physics B}

\begin{document}

\begin{frontmatter}



\title{Mapping-based Hard-constrained Physics-Informed Neural Networks for unbounded wave problems} %% Article title


\author[a]{Tao Zhang}
\author[a,b]{Hanshu Chen}
\author[c]{Ilia Marchevsky}
\author[a]{Zhuojia Fu\corref{cor1}}

\cortext[cor1]{Corresponding author.}

\address[a]{Center for Numerical Simulation Software in Engineering and Sciences, College of Mechanics and Engineering Science, Hohai University, Nanjing, Jiangsu, 211100, China}
\address[b]{State Key Laboratory of Structural Analysis for Industrial Equipment, Optimization and CAE Software, Department of Engineering Mechanics, Dalian University of Technology, Dalian 116024, China}
\address[c]{Bauman Moscow State Technical University, 2-nd Baumanskaya st., 5, 105005, Moscow, Russia}

%% Abstract
\begin{abstract}
The aim of this paper is to introduce a Mapping-based Hard-constrained Physics-Informed Neural Network (MH-PINN) for efficiently and accurately solving unbounded wave problems. First, we propose a coordinate mapping technique that compactifies the infinite physical domain into a finite computational space. This effectively resolves the sampling difficulties inherent to standard PINNs in unbounded regions. Additionally, it avoids the artificial truncation errors introduced by traditional methods such as perfectly matched layers. Second, we design a physics-based hard-constrained network structure that automatically satisfies both the inner boundary conditions and the far-field radiation conditions. This structure eliminates boundary loss terms, yielding high computational efficiency and fast convergence, which effectively addresses the challenges of high-frequency problems. Third, we introduce an inverse factor correction for boundary coefficients to address the influence of asymptotic factors,which makes the method highly geometrically adaptable. Finally, we present numerical examples covering various acoustic radiation and scattering scenarios as well as elastic dynamics scenarios to demonstrate the efficiency and accuracy of our algorithm.It highlights its potential for broader applications in the field of computational wave dynamics.
\end{abstract}

%% Keywords

\end{frontmatter}



\section{Introduction}
\label{sec:introduction}

Solving partial differential equations in unbounded domains remains a fundamental challenge in scientific computing, particularly for acoustic radiation and scattering problems where physical fields extend to infinity \cite{li2012analysis,chandler2012numerical,alves2024wave}. Accurately capturing far-field asymptotic behaviors is crucial. Traditional grid-based methods, such as the finite element method \cite{zienkiewicz1977finite,jiang2013numerical,li2015hybrid,vanderEerden1996}, finite difference method \cite{zhang2023finite,thomas2013numerical}, and meshfree methods \cite{rabczuk2019extended}, typically truncate the domain using artificial absorbing boundary conditions \cite{engquist1977absorbing} or perfectly matched layers \cite{berenger1994perfectly}. However, these treatments often introduce spurious numerical reflections and face the ``curse of dimensionality'' \cite{hao2022physics} when resolving high-frequency oscillations over vast domains. Alternative approaches, like the boundary element method  \cite{wu1995direct,chandler2004high} and spectral methods \cite{xiu2007efficient,cho2019spectrally,shen2009some}, struggle with dense matrix operations for heterogeneous media \cite{liu2011recent} and lack geometric adaptability for complex scatterers, respectively.

Physics-Informed Neural Networks (PINNs)\cite{raissi2019physics,karniadakis2021physics,jeong2025advanced,zhongkai2024pinnacle,luo2025physics,tang2022extrinsic,fu2024physics} have emerged as a promising mesh-free solver, yet their direct application to infinite acoustic domains presents severe bottlenecks, despite some pioneering domain-extension efforts \cite{xia2023spectrally,ren2024seismicnet}. First, generating effective training collocation points across an infinite space is practically infeasible, leading to poor generalization in the far field. Second, the ``spectral bias'' of fully connected networks hinders the learning of high-frequency acoustic wave features \cite{krishnapriyan2021characterizing,rahaman2019spectral,wang2021eigenvector}. Finally, current PINNs predominantly enforce radiation and boundary conditions via soft penalty terms. Because far-field values are orders of magnitude smaller than near-field residuals, this soft constraint mechanism exacerbates ``gradient competition'' \cite{wang2021understanding} and often traps the optimization in local minima.

To address these challenges, we develop a Mapping-based Hard-constrained Physics-Informed Neural Network (MH-PINN) framework. By integrating spectral mapping principles with the adaptive nature of deep learning, the proposed MH-PINN offers a robust approach for unbounded problems. Our framework features two primary innovations:
\begin{itemize}
    \item \textbf{Domain Compactification:} A coordinate transformation maps the infinite physical domain into a bounded parametric space. This completely resolves the unbounded sampling dilemma, allowing a finite set of collocation points to cover the entire infinite space efficiently.
    \item \textbf{Physical Factorization for Hard Constraints:} We design a constructive solution hypothesis that strictly couples a geometric distance function  with an analytical attenuation factor. By ``hard-coding'' these constraints directly into the network architecture, we entirely eliminate the competing boundary loss terms, simplifying the high-dimensional non-convex optimization landscape to depend solely on the governing PDE residuals.
\end{itemize}

The remainder of this paper is organized as follows: Section \ref{sec:preliminaries} formulates the governing equation for the unbounded acoustic problem and analyzes the inherent limitations of the standard PINN framework. Section \ref{sec:mh_pinn} introduces the proposed MH-PINN framework, detailing the analytical coordinate mapping and the physical hard-constraint architecture. Section \ref{sec:results} systematically validates the method's accuracy, stability, and geometric adaptability through acoustic radiation and scattering benchmarks, including scenarios with heterogeneous media and non-circular boundaries. Section \ref{sec:limitations} discusses current limitations, and Section \ref{sec:conclusions} concludes the paper with future perspectives.
\section{Problem Statement and Limitations of Standard PINNs}
\label{sec:preliminaries}

\subsection{Problem Statement}

We consider the Helmholtz equation, which is fundamental to acoustic radiation, scattering and elastodynamics, defined on an unbounded exterior domain $\Omega \subset \mathbb{R}^d$\begin{equation}
    \nabla^2 u(\mathbf{x}) + k^2(\mathbf{x}) u(\mathbf{x}) = f(\mathbf{x}), \quad \mathbf{x} \in \Omega
    \label{eq:helmholtz}
\end{equation}
\noindent where $\nabla^2$ is the Laplacian operator, $k$ is the wavenumber, and $f(\mathbf{x})$ is the source term. The complex-valued field $u(\mathbf{x})$ must satisfy specific physical boundary conditions on the internal boundary $\Gamma = \partial\Omega$. In acoustic radiation and sound-soft scattering problems, this is typically a Dirichlet condition:
\begin{equation}
    u(\mathbf{x}) = g_D(\mathbf{x}), \quad \mathbf{x} \in \Gamma
    \label{eq:dirichlet}
\end{equation}
\noindent where $g_D(\mathbf{x})$ represents the prescribed boundary values. Conversely, for sound-hard acoustic scatterers or traction-free boundaries in elastodynamics (e.g., SH wave scattering), the field is governed by a Neumann boundary condition:
\begin{equation}
    \frac{\partial u(\mathbf{x})}{\partial \mathbf{n}} = g_N(\mathbf{x}), \quad \mathbf{x} \in \Gamma
    \label{eq:neumann}
\end{equation}
\noindent where $\mathbf{n}$ is the outward unit normal vector to the boundary $\Gamma$, and $g_N(\mathbf{x})$ denotes the prescribed normal derivative or traction forces. Furthermore, to ensure a unique and physically meaningful solution in the infinite far field, $u(\mathbf{x})$ must satisfy the Sommerfeld radiation condition:
\begin{equation}
    \lim_{r \to \infty} r^{1/2} \left( \frac{\partial u}{\partial r} - iku \right) = 0, \quad r = |\mathbf{x}|
\label{eq:sommerfeld}
\end{equation}
\noindent where $r$ denotes the radial distance from the origin and $i$ is the imaginary unit. This condition physically implies that acoustic or elastodynamic scattered waves must propagate outward toward infinity without reflecting back from the boundary of the domain.

\subsection{Solving the Helmholtz Equation using Standard PINNs}
Physics-Informed Neural Networks (PINNs) have emerged as a promising mesh-free approach for solving partial differential equations\cite{raissi2019physics}. In the standard PINN framework, a fully connected neural network $\hat{u}(\mathbf{x};\theta)$ parameterized by weights $\theta$ is employed to approximate the latent acoustic field. Instead of relying on traditional mesh discretization, PINNs substitute the network output into the governing PDE using automatic differentiation\cite{raissi2019physics}.

To solve the unbounded Helmholtz problem, standard PINNs formulate the training process as an optimization problem that minimizes a composite soft-constrained loss function\cite{nair2024acoustic}:
\begin{equation}
    \mathcal{L}_{total}=\lambda_{PDE}\mathcal{L}_{PDE}+\lambda_{BC}\mathcal{L}_{BC}+\lambda_{RAD}\mathcal{L}_{RAD}
    \label{eq:loss_pinn}
\end{equation}
where $\mathcal{L}_{PDE}$ represents the mean squared error of the Helmholtz equation residual evaluated at interior collocation points, $\mathcal{L}_{BC}$ enforces the Dirichlet boundary condition on $\Gamma$, and $\mathcal{L}_{RAD}$ penalizes the violation of the radiation condition (or an absorbing boundary condition) at a truncated far-field boundary. The hyperparameters $\lambda_{PDE}$, $\lambda_{BC}$, and $\lambda_{RAD}$ are penalty weights used to balance the competing loss terms.

\subsection{Challenges of Standard PINNs}
While standard PINNs are effective for many bounded PDE problems, their direct application to Eqs. (\ref{eq:helmholtz})--(\ref{eq:sommerfeld}) faces severe computational bottlenecks:

\begin{itemize}
    \item \textbf{Unbounded domain sampling:} Generating effective training collocation points uniformly across an infinite space is practically infeasible. Imposing artificial truncation boundaries requires massive sampling in vast domains, leading to poor generalization in the far field.
    \item \textbf{Spectral bias:} The implicit "spectral bias" of fully connected networks impedes the learning of high-frequency oscillatory wave features.as noted in\cite{wang2022and},The network tends to fit low-frequency components first, struggling to capture the rapid phase variations of acoustic waves at higher wavenumbers.
    \item \textbf{Gradient competition:} Enforcing radiation and inner boundary conditions via soft penalty terms creates a highly non-convex optimization landscape. Because far-field values are often orders of magnitude smaller than near-field residuals, this mechanism exacerbates "gradient competition," trapping the optimizer in local minima and severely hindering convergence.
\end{itemize}

These inherent limitations motivate the development of MH-PINN proposed in this work, which integrates coordinate compactification and physical hard constraints to resolve the above issues fundamentally.

\section{The Proposed MH-PINN Framework}
\label{sec:mh_pinn}

In this section, we present the full formulation and network design of the MH-PINN. Starting from a generalized coordinate transformation for unbounded domain compactification, we then construct a physics-based hard-constrained solution ansatz that automatically satisfies arbitrary inner boundary conditions and far-field radiation conditions by structure. Finally, we introduce the dual-stream network architecture and a purely PDE-based loss function, yielding a concise, stable, and high-accuracy solver for unbounded wave problems.

\subsection{Coordinate Mapping for Domain Compactification}
\label{subsec:mapping}

In order to solve the computation and sampling difficulties caused by unbounded
 physical regions $\Omega = \{(r, \theta) \mid r \in [R_{in}, \infty), \theta \in [0, 2\pi)\}$, 
 we introduce a coordinate mapping strategy to convert the infinite external region into a finite computational domain $\tilde{\Omega}$. 
 This method avoids the truncation error inherent in artificial boundary conditions and
  allows neural networks to be sampled and trained on a compact domain. 
  Shen et al. \cite{shen2009some} pointed out that a suitable mapping can 
  transform differential operators on unbounded domains into singular or non-singular operators on bounded domains.
   Considering that waves usually exhibit algebraic decay at infinity, and to avoid the numerical overflow problem that exponential mapping may cause in the far field, this paper adopts a robust algebraic mapping strategy.

\subsubsection{Mapping Function Formula}
Let $\mathbf{x} = (r, \theta)$ represent the coordinates in the physical domain and $\zeta = (\xi, \theta)$ represent the coordinates in the computational domain $\tilde{\Omega}$. An algebraic transformation is used to compactify the semi-infinite interval in the physical domain $[R_{in}, \infty)$ to a finite interval in the computational domain $[-1, 1)$:
\begin{equation}
    r(\xi) = R_{in} + L\frac{1+\xi}{1-\xi}, \quad \xi \in [-1, 1)
\end{equation}
where $R_{in}$ is a reference radius forming a background mapping, and $L > 0$ is a scaling parameter that controls the distribution of collocation points. In this mapping, $r = R_{in}$ when $\xi = -1$, and $r \to \infty$ when $\xi \to 1$.

\begin{figure}[htbp]
    \centering
    \includegraphics[width=\linewidth]{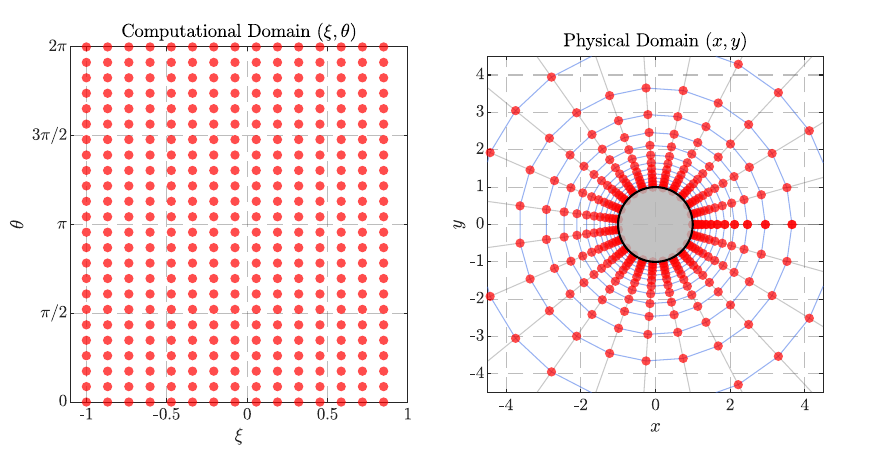}
    \caption{Visualization of spatial mapping transformation from computational domain $(\xi, \theta)$ to physical domain $(x, y)$.}
    \label{fig:spatial_mapping}
\end{figure}

It is crucial to emphasize that in our generalized framework, $R_{in}$ serves merely as a background mapping radius and does not need to conform to the actual complex boundary of the scatterer. This cleanly decouples the spatial compactification from the geometric representation. As shown in Figure \ref{fig:spatial_mapping}, this algebraic mapping compresses the semi-infinite region into a finite interval, achieving a spatial discretization effect of "near-field density and far-field sparsity." This strategy maximizes computational efficiency while entirely bypassing the spurious reflection problem caused by artificial truncation.

\subsubsection{Transformation of Differential Operators and Governing Equations}
To train the neural network in a bounded domain $\tilde{\Omega}=[-1,1) \times [0, 2\pi)$, the differential operators in the original physical problem must be reconstructed. By applying the chain rule, the radial derivatives can be expressed as functions of $\xi$:
\begin{equation}
    \frac{\partial u}{\partial r} = \frac{1}{J(\xi)} \frac{\partial u}{\partial \xi}, \quad 
    \frac{\partial^2 u}{\partial r^2} = \frac{1}{J(\xi)^2} \frac{\partial^2 u}{\partial \xi^2} - \frac{J'(\xi)}{J(\xi)^3} \frac{\partial u}{\partial \xi}
\end{equation}
where $J(\xi) = \frac{dr}{d\xi} = \frac{2L}{(1-\xi)^2}$ and $J'(\xi) = \frac{4L}{(1-\xi)^3}$. The Jacobian determinant $J(\xi)$ encodes the geometric stretching information, acting as a mathematical bridge connecting the two spaces.

Substituting these transformations into the standard Helmholtz equation yields the modified governing equation defined on the computational domain $(\xi, \theta)$:
\begin{equation}
    \frac{1}{J^2} \frac{\partial^2 u}{\partial \xi^2} + \left( \frac{1}{r(\xi)J} - \frac{J'}{J^3} \right) \frac{\partial u}{\partial \xi} + \frac{1}{r(\xi)^2} \frac{\partial^2 u}{\partial \theta^2} + k^2 u = 0
    \label{eq:transformed_pde}
\end{equation}

Within this framework, the core training objective is transformed into minimizing the residual of Eq.(\ref{eq:transformed_pde}) within the bounded computational domain.

\subsection{Physical Hard Constraints Construction}
\label{subsec:hard_constraints}

When dealing with acoustic and elastodynamic wave problems in unbounded regions, numerical solutions must simultaneously satisfy the specific physical conditions at the inner boundary and the Sommerfeld radiation conditions at infinity. To avoid gradient competition caused by soft-constraint penalties, we construct a generalized, geometry-agnostic hard-constraint framework based on physical factor decomposition. We explicitly decompose the predicted solution $\hat{u}(\mathbf{x})$ into a product of a rapidly varying asymptotic factor $\Phi(\mathbf{x})$ and a slowly varying neural envelope $\mathcal{E}(\mathbf{x})$:

\begin{equation}
    \hat{u}(\mathbf{x}) = \underbrace{\Phi(\mathbf{x})}_{\text{Far-field Asymptotic}} \cdot \underbrace{\mathcal{E}(\mathbf{x})}_{\text{Neural Envelope}}
\end{equation}

\subsubsection{Far-field Asymptotic Factor for Radiation Condition}
In unbounded domain problems, the physical domain's $|\mathbf{x}| \to \infty$ typically corresponds to the mapped compact computational coordinate $\xi \to 1$. Because the mapping Jacobian exhibits a singularity at the far-field limit, directly fitting the neural network and applying the derivative constraint of the radiation condition leads to numerical instability. 

To circumvent this singularity and satisfy the radiation condition, we embed the physical prior knowledge directly into $\Phi(\mathbf{x})$. For outward-propagating waves in a $d$-dimensional space, we define:
\begin{equation}
    \Phi(\mathbf{x}) = \frac{e^{ik(|\mathbf{x}| - R_{in})}}{|\mathbf{x}|^{(d-1)/2}}
\end{equation}
where $R_{in}$ is the background reference radius defined in the mapping.
 For a two-dimensional problem ($d=2$), noting that $r = |\mathbf{x}|$,
  this simplifies to $\Phi(r) = r^{-1/2} e^{ik(r - R_{in})}$ as depicted in our
   network architecture. This factor analytically describes the geometric attenuation and phase oscillation in the far field, automatically and strictly satisfying the Sommerfeld radiation condition mathematically, regardless of the network's initialization.

\subsubsection{Exact Distance Function and Its Geometric Properties}
To accommodate arbitrary and complex inner boundaries $\Gamma$ independently of the background mapping, we introduce the theory of exact distance functions. We construct a normalized implicit distance function $d(\mathbf{x})$ for the specific target geometry. 

This distance function strictly satisfies two elegant and crucial geometric properties:
(1) Zero-value property: $d(\mathbf{x}) = 0$ for all points $\mathbf{x} \in \Gamma$;
(2) Normal consistency: $\nabla d(\mathbf{x}) = \mathbf{n}(\mathbf{x})$ for all points $\mathbf{x} \in \Gamma$, where $\mathbf{n}$ is the outward unit normal vector at that point.

By utilizing Boolean operations, the exact distance function $d(\mathbf{x})$ can be analytically constructed for any combination of complex geometries, entirely decoupling the representation of specific boundary shapes from the neural network optimization process and the spatial mapping.

\subsubsection{Hard Constraint Construction for Dirichlet and Neumann Boundaries}
Based on the geometric properties provided by the exact distance function $d(\mathbf{x})$, we can construct penalty-free, exact hard constraints for various physical boundary conditions. Here, $\mathcal{N}(\xi, \theta; \mathbf{w})$ represents the neural network output parameterized by weights $\mathbf{w}$, which takes the mapped compact coordinates as inputs to ensure expressivity within the infinite domain.

For Dirichlet boundary conditions (e.g., acoustic radiation or sound-soft scattering problems), the field must satisfy $u(\mathbf{x}) = g_D(\mathbf{x})$. Leveraging the zero-value property of $d(\mathbf{x})$, we parameterize the neural envelope $\mathcal{E}(\mathbf{x})$ as:
\begin{equation}
    \mathcal{E}(\mathbf{x}) = \frac{g_D(\mathbf{x})}{\Phi(\mathbf{x})} + d(\mathbf{x}) \cdot \mathcal{N}(\xi, \theta; \mathbf{w})
\end{equation}

Because $d(\mathbf{x}) \equiv 0$ on the boundary $\Gamma$, the second term on the right side naturally vanishes, inherently yielding $\hat{u}(\mathbf{x})\big|_\Gamma = g_D(\mathbf{x})$ for the total field. This ansatz strictly and automatically imposes the Dirichlet condition on arbitrary geometries.

For Neumann boundary conditions (e.g., sound-hard scatterers or traction-free boundaries), the normal derivative of the field is constrained, i.e., $\frac{\partial u}{\partial n} = \nabla u \cdot \mathbf{n} = g_N(\mathbf{x})$. The normal derivative of the total field $\hat{u} = \Phi \mathcal{E}$ on the boundary expands to:
\begin{equation}
    \frac{\partial \hat{u}}{\partial n}\bigg|_\Gamma = \Phi (\nabla \mathcal{E} \cdot \mathbf{n}) + \mathcal{E} (\nabla \Phi \cdot \mathbf{n}) = g_N(\mathbf{x})
\end{equation}

To satisfy this condition implicitly, the envelope function $\mathcal{E}$ must possess a specific target normal derivative $h(\mathbf{x})$ on the boundary:
\begin{equation}
    h(\mathbf{x}) = \frac{g_N(\mathbf{x}) - \mathcal{N}(\xi, \theta; \mathbf{w}) \left[ \nabla \Phi(\mathbf{x}) \cdot \nabla d(\mathbf{x}) \right]}{\Phi(\mathbf{x})}
\end{equation}

Utilizing the normal consistency property $\nabla d = \mathbf{n}$ on the boundary, we propose an analytical Taylor-expansion compensation method to construct a first-order shielded envelope:
\begin{equation}
    \mathcal{E}(\mathbf{x}) = \mathcal{N}(\xi, \theta; \mathbf{w}) - d(\mathbf{x}) \Big[ \nabla \mathcal{N}(\xi, \theta; \mathbf{w}) \cdot \nabla d(\mathbf{x}) - h(\mathbf{x}) \Big]
\end{equation}

By embedding this shielded envelope directly into the forward pass, the framework mathematically and intrinsically guarantees the exact satisfaction of Neumann boundary conditions for completely arbitrary shapes. It should be emphasized that the gradient operator $\nabla$ in the above formulation denotes the spatial derivative with respect to the physical coordinates $\mathbf{x}$. The term $\nabla \mathcal{N}(\xi, \theta; \mathbf{w})$ implicitly utilizes the chain rule through the inverse mapping functions, which is seamlessly handled by the automatic differentiation machinery of deep learning frameworks.

\subsection{Network Architecture and No-Boundary Loss Function}
\begin{figure}[H]
    \centering
    \includegraphics[width=\linewidth]{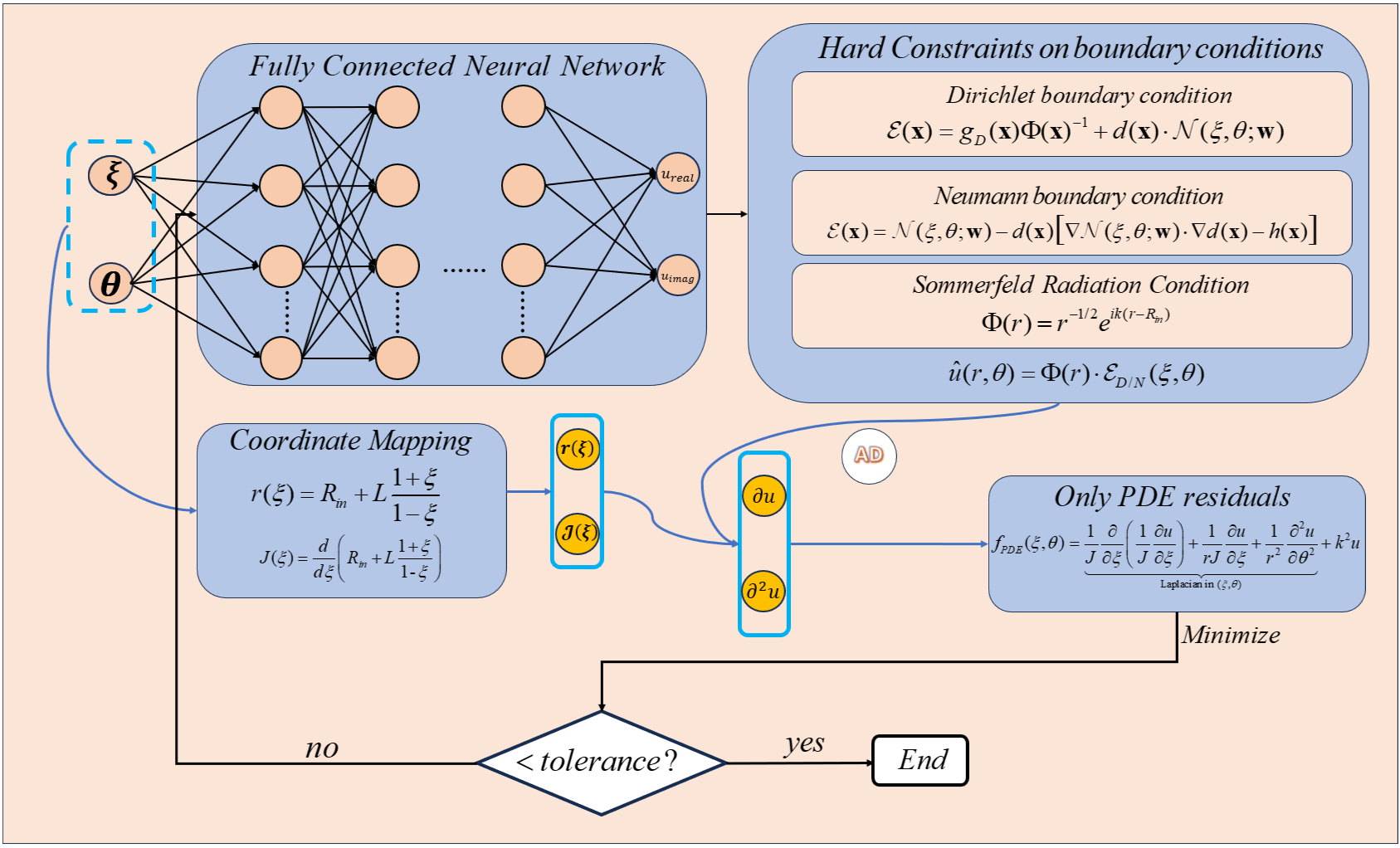} 
    \caption{Framework diagram of the proposed mapping-based hard-constrained PINN. (Note: $R_{in}$ and the 1D analytical PDE formulation in the figure represent the baseline radial scenario, which is universally generalized via Automatic Differentiation and $r_b(\theta)$ in the code implementation).}
    \label{fig:network_architecture}
\end{figure}

As illustrated in Figure~\ref{fig:network_architecture}, the generalized MH-PINN employs a tightly coupled "dual-stream" architecture. The Fully Connected Neural Network stream takes the mapped compact coordinates $(\xi, \theta)$ as inputs to ensure sufficient expressivity across the infinite domain, outputting the raw network variables $\mathcal{N}(\xi, \theta; \mathbf{w})$ (including real and imaginary parts). Simultaneously, the Coordinate Mapping stream analytically calculates $r(\xi)$ and the Jacobian $J(\xi)$, providing geometric scaling factors without participating in network gradient updates.

In the subsequent "Hard Constraints" module, the framework fuses the network outputs $\mathcal{N}$, the far-field asymptotic factor $\Phi(\mathbf{x})$, and the exact geometric distance function $d(\mathbf{x})$ (along with its gradient $\nabla d(\mathbf{x})$) based on the rigorous mathematical ansätze derived in Section 3.2. This fusion constructs the final predicted field $\hat{u}(r, \theta)$, universally satisfying either Dirichlet or Neumann boundary conditions and the Sommerfeld radiation condition by design.

Thanks to this inherently physics-compliant structural design, all boundary data terms and penalty hyperparameters are completely eliminated. The total loss function simplifies to include only the mean squared error of the PDE residuals evaluated in the mapped computational domain:
\begin{equation}
    \mathcal{L}_{total}(\mathbf{w}) = \mathcal{L}_{PDE} = \frac{1}{N} \sum_{i=1}^{N} \left| f_{PDE}(\xi_i, \theta_i; \hat{u}) \right|^2
\end{equation}
where $f_{PDE}(\xi, \theta)$ is the analytically mapped Helmholtz equation residual (as explicitly shown in Figure~\ref{fig:network_architecture}). This elegant unconstrained minimization allows the optimizer to bypass gradient competition, achieving stable and highly accurate unsupervised learning across the entire unbounded domain.

\section{Numerical Examples}
\label{sec:results}
\subsection{Acoustic Radiation Problem}

To comprehensively evaluate the performance of MH-PINN in solving unbounded acoustic radiation, we designed three distinct sets of numerical experiments. The first set involves a wideband robustness test, examining the algorithm's stability as the wavenumber $k$ varies from $1\mathrm{\ m}^{-1}$ to $10\mathrm{\ m}^{-1}$. The second set focuses on a high-frequency limit test at $k = 20\mathrm{\ m}^{-1}$ to assess performance in strongly oscillatory scenarios. The third set considers a more complex case with a spatially varying wavenumber, $k(\mathbf{x})$, to simulate acoustic radiation in inhomogeneous media, further testing the method's adaptability to non-uniform physical properties.

We first consider the two-dimensional acoustic radiation from an infinitely long cylinder. The governing Helmholtz equation and the corresponding boundary conditions follow the definitions provided in Section 2, defined over the domain $r \in [1, \infty)$ and $\theta \in [0, 2\pi)$. The inner boundary is defined as the unit circle $r=1$, where a constant Dirichlet boundary condition $u(1, \theta) =100\mathrm{~Pa}$ is imposed. At the boundary at infinity, the field satisfies the Sommerfeld radiation condition, ensuring the outward propagation of sound waves without reflection.For such an axisymmetric radiation problem, the exact analytical solution is well-established \cite{morse1986theoretical} and expressed in terms of the Hankel function of the first kind and zero order:
\begin{equation}
    u_{exact}(r) = u_0 \frac{H_0^{(1)}(kr)}{H_0^{(1)}(k)}
\end{equation}
\noindent This analytical solution can be used as the true value to verify the accuracy of the numerical results.

\begin{figure}[H]
    \centering
    \begin{minipage}{0.85\textwidth}
        \centering
        \includegraphics[width=\textwidth]{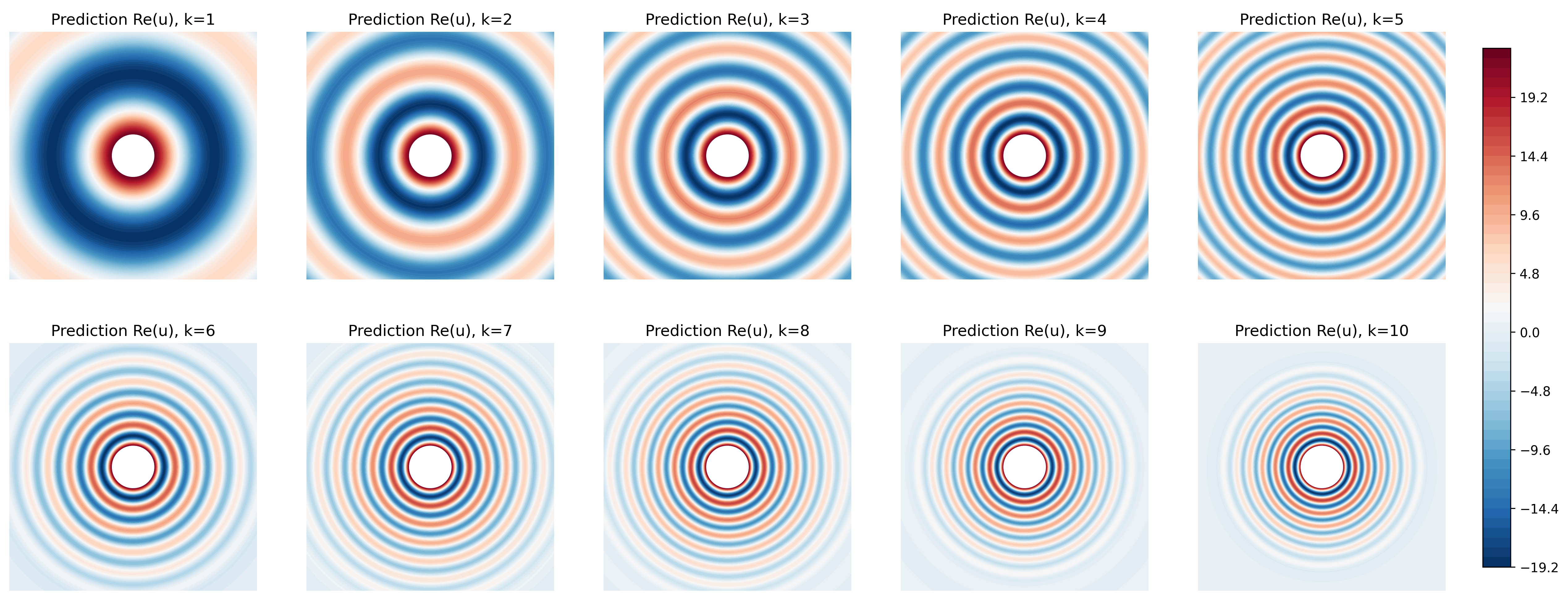}
        \caption{Predicted real part of sound pressure by Physics-Informed Neural Network}
        \label{fig:pred_pinn}
    \end{minipage}
    
    \vspace{10pt} % 调整这个值来控制间距
    
    \begin{minipage}{0.85\textwidth}
        \centering
        \includegraphics[width=\textwidth]{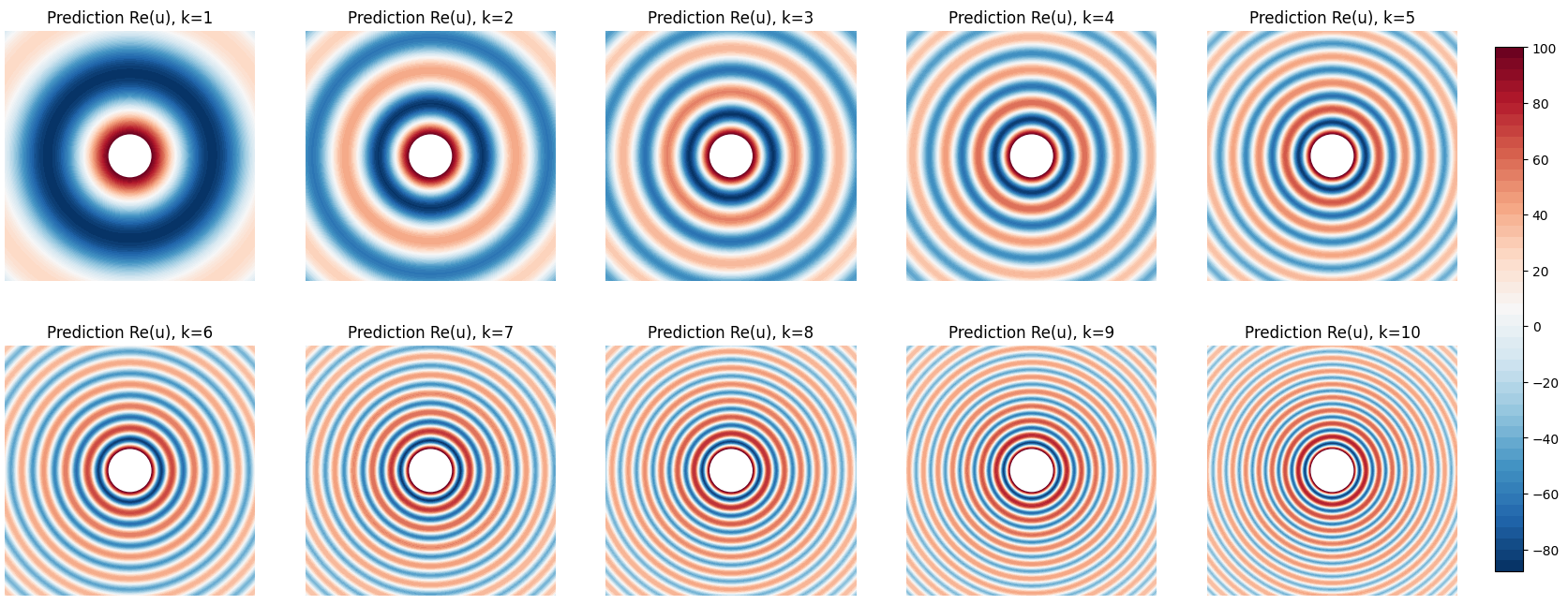}
        \caption{Predicted real part of sound pressure by mapping-based hard-constrained Physics-Informed Neural Network}
        \label{fig:pred_mh-pinn}
    \end{minipage}
\end{figure}

Figures~\ref{fig:pred_pinn} and \ref{fig:pred_mh-pinn} respectively display the real-part contour plots of the sound pressure field predicted by the standard PINN and MH-PINN. Observation of Figure 4 reveals that under low wave number conditions ($k \leq 5$), the standard PINN effectively captures acoustic wave characteristics. However, as the wave number $k$ increases, the standard PINN's predictions rapidly deteriorate, with acoustic field stripes becoming blurred. This indicates that the standard PINN suffers from severe spectral bias when handling high-frequency oscillatory problems, struggling to converge towards the correct solution. In contrast, Figure 5 demonstrates that MH-PINN clearly and accurately reconstructs the concentric circular acoustic radiation pattern across the entire test frequency range ($k = 1$ to $10$). Even at $k=10$, the peaks and troughs remain distinctly discernible without numerical dissipation or distortion.

\begin{figure}[htbp] % 如果你之前用了 float 宏包，这里也可以换成 [H] 
    \centering

    % 左侧子图 (a)：误差分析
    \begin{subfigure}[b]{0.48\textwidth}
        \centering
        \includegraphics[width=\textwidth]{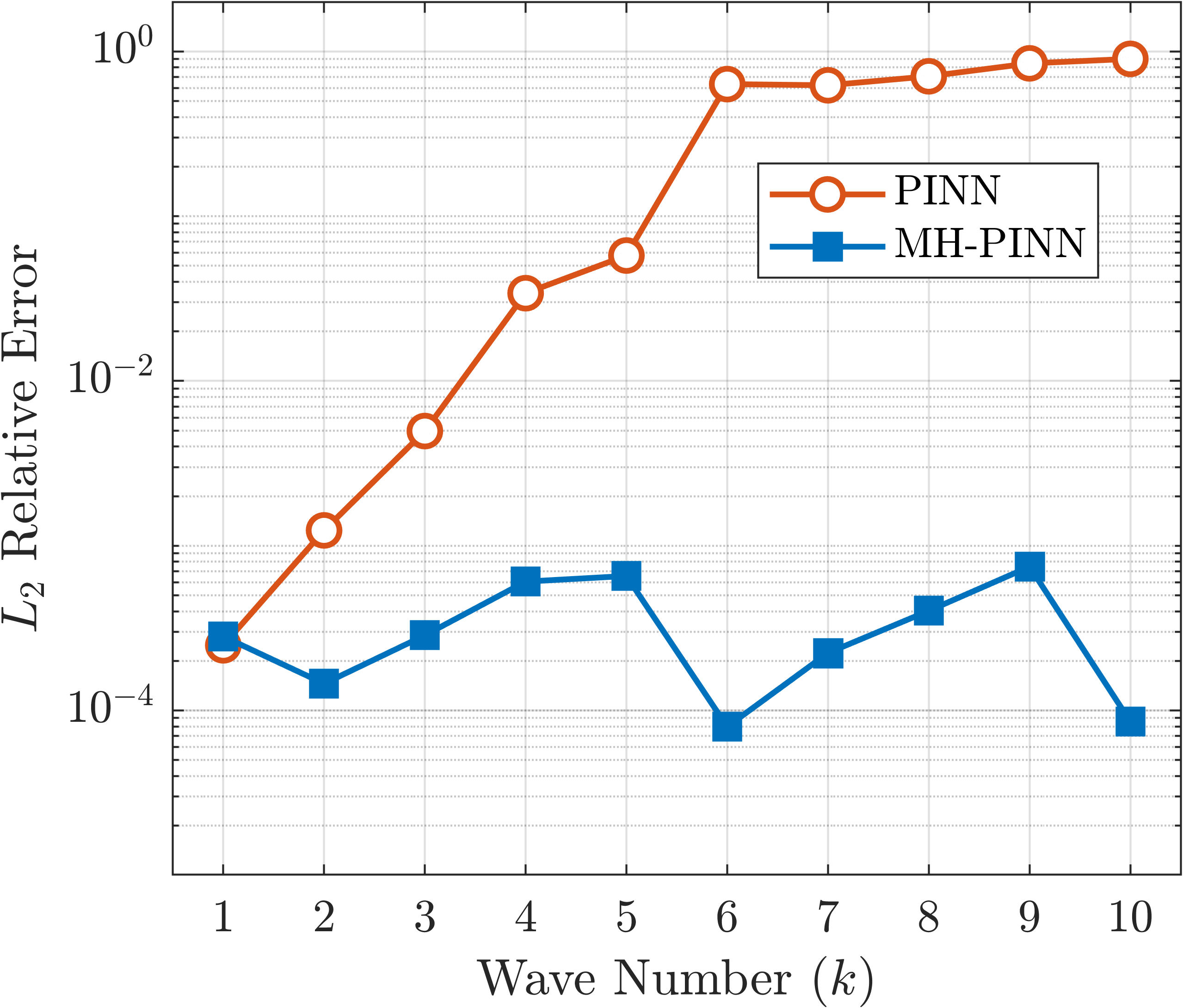}
        \caption{Computational Error Analysis}
        \label{fig:error}
    \end{subfigure}
    \hfill % hfill 会自动在两个子图之间插入适当的空白间距 
    % 右侧子图 (b)：效率分析
    \begin{subfigure}[b]{0.48\textwidth}
        \centering
        \includegraphics[width=\textwidth]{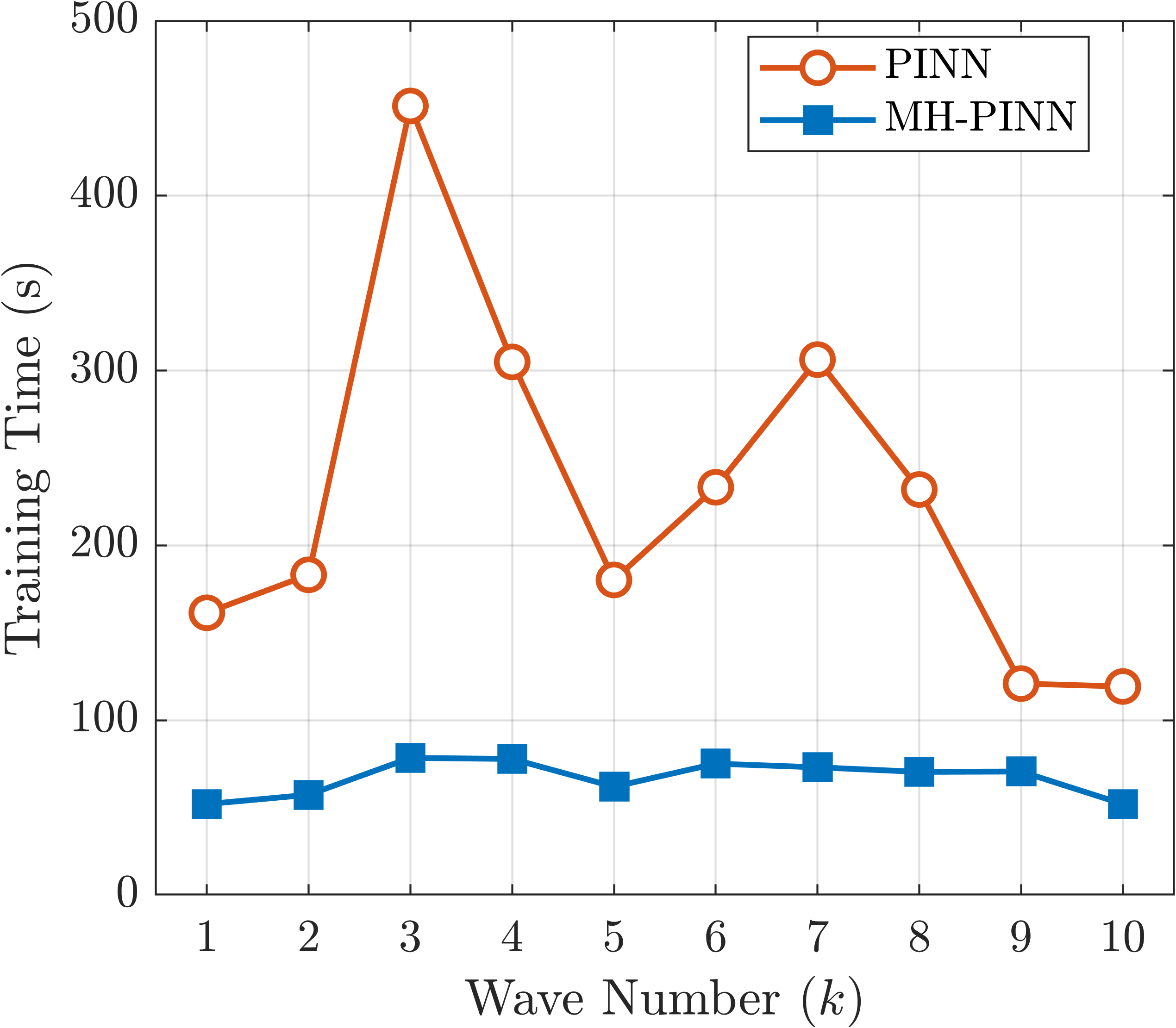}
        \caption{Computational Efficiency Analysis}
        \label{fig:efficiency}
    \end{subfigure}
    
% 整个大图的图注
    \caption{Comparison of computational performance between standard PINN and MH-PINN across different wavenumbers: (a) Error analysis; (b) Efficiency analysis.}
    \label{fig:radiation_comparison}
\end{figure}

Figure \ref{fig:radiation_comparison} visually quantifies the relative errors and computational time consumption of both methods across varying wavenumbers from $k = 1$ to $10$.As shown in Figure 6(a), the PINN’s error curve displays a significant upward trend as the wave number increases. By $k = 6$, the PINN's error skyrockets to $6.33 \times 10^{-1}$, indicating a complete failure of the model, and eventually deteriorates to $9.01 \times 10^{-1}$ at $k = 10$. Conversely, the MH-PINN error curve demonstrates exceptional numerical stability. Benefiting from the coordinate mapping and hard-constraint mechanisms, its $L_2$ relative error fluctuates within an extremely narrow band, consistently remaining between $7.97 \times 10^{-5}$ and $7.46 \times 10^{-4}$ across all tested wavenumbers. This stark contrast highlights the algorithm's strong robustness against frequency variations, showing no divergence trend even in high-frequency scenarios.Regarding computational efficiency, Figure 6(b) reveals that the training time for MH-PINN is significantly lower and more stable than that of the PINN. The PINN suffers from optimization bottlenecks, with training times fluctuating violently and peaking at over 450 seconds (at $k = 3$). In stark contrast, MH-PINN maintains a highly efficient training profile, requiring only 52.1 to 78.5 seconds across all scenarios. This shows that automatically satisfying boundary and radiation conditions through hard constraints significantly reduces the optimization difficulty of the loss function, enabling the network to bypass gradient competition, find the global optimum faster, and maintain higher accuracy.

\begin{figure}[htbp]
    \centering
    \includegraphics[width=0.85\textwidth]{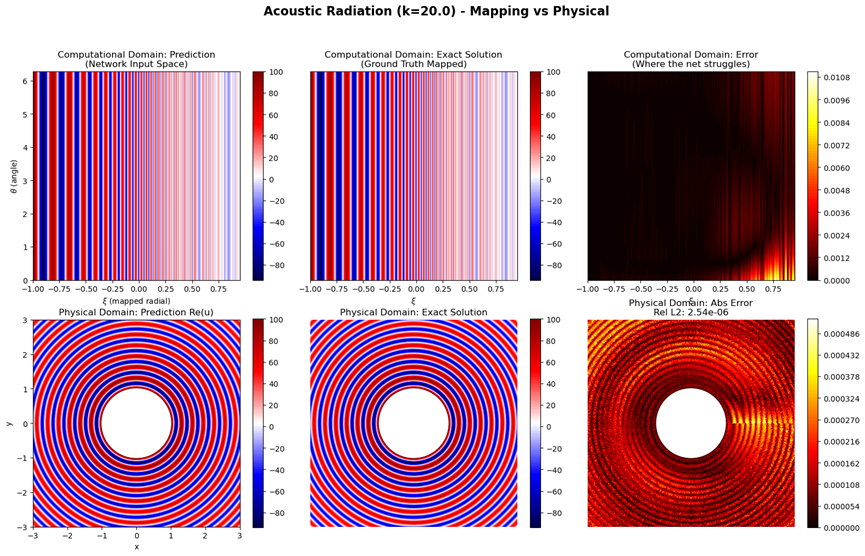}
    \caption{Predicted solution, exact solution, and error for the MH-PINN global computational domain and test physical domain at wave number k=20}
    \label{fig:limiting performance}
\end{figure}

To further investigate the limiting performance of MH-PINN at higher frequencies, we simulated the strongly oscillatory case of $k=20$.This constitutes an exceptionally challenging example, where conventional meshfree methods typically require extremely high sampling densities. Figure \ref{fig:limiting performance} presents detailed results for MH-PINN across the global computational domain and the test physical domain.

The left and middle columns respectively compare the network-predicted solution with the analytical solution. It can be observed that after compressing the unbounded physical domain $r \in [1, \infty)$ into a finite computational domain $\xi \in [-1, 1]$ via coordinate mapping, the oscillatory solution---which originally decayed with radius---is transformed into regular periodic stripes within the computational domain. MH-PINN accurately captures this high-frequency oscillatory structure.

The right column displays the absolute error distribution across the global computational domain and the test physical domain. The final quantitative assessment reveals that at high wave number $k = 20$, the MH-PINN achieves a relative error of $2.54 \times 10^{-6}$. This numerical experimental result demonstrates that the mapping-based hard constraint strategy not only successfully resolves the unbounded domain truncation problem but also effectively overcomes numerical dispersion and contamination effects in solving high-frequency Helmholtz equations.

In the third example, we demonstrate the proposed method's capability to address acoustic radiation problems in non-uniform media. Unlike the constant wave number case, wave propagation in non-uniform media poses significant computational challenges due to spatially varying wavelengths and the accumulation of phase errors over long distances.

We consider a two-dimensional, unbounded acoustic radiation problem with spatially varying coefficients, governed by the two-dimensional Helmholtz equation:
\begin{equation}
    \nabla^2 u(r, \theta) + k^2(r)u(r, \theta) = 0, \quad r \in [1, \infty)
    \label{eq:helmholtz_varying}
\end{equation}

where $u$ denotes the complex sound pressure. The wave number distribution is defined as:
\begin{equation}
    k(r) = k_\infty \left( 1 + \alpha e^{-(r-1)/d} \right)
    \label{eq:wavenumber_distribution}
\end{equation}

where $k$ denotes the far-field background wave number, $\alpha$ denotes the perturbation intensity coefficient, and $d$ denotes the decay scale factor. This problem constitutes a challenging multiscale scenario, implying high wave numbers and short wavelengths near the cylindrical surface where the medium is relatively `dense', with its properties decaying exponentially towards the uniform background state with increasing distance.

The internal boundary conditions are identical to those in the preceding section. To address the issue of varying wave numbers, we introduce the WKB approximation principle as a hard constraint on our far-field radiation boundary. Rather than employing a linear phase term, we adopt an integral phase $\Phi(r)$ to describe the spatially varying refractive index:
\begin{equation}
    \Phi(r) = \int_1^r k(\tau) d\tau = k_\infty(r-1) + k_\infty \alpha d \left( 1 - e^{-(r-1)/d} \right)
    \label{eq:integral_phase}
\end{equation}

Therefore, the solution is approximated by the network output $\mathcal{N}$ constructed using the proposed coordinate mapping and hard constraints:
\begin{equation}
    \hat{u}(r, \theta) = \frac{e^{i\Phi(r)}}{\sqrt{r}} \left[ u_{bc} + (\xi(r) + 1)\mathcal{N}(\xi, \theta; \mathbf{w}) \right]
    \label{eq:network_output}
\end{equation}

We employed Latin hypercube sampling within the computational domain $\Omega_{comp} = [-1, 1] \times [0, 2\pi]$ to generate 15,000 data points. In the absence of an analytical solution, a reference solution was generated using a high-order finite difference method over a larger truncated domain to serve as the ground truth for verifying the method's predictive accuracy.

\begin{figure}[htbp]
    \centering
    \includegraphics[width=0.85\textwidth]{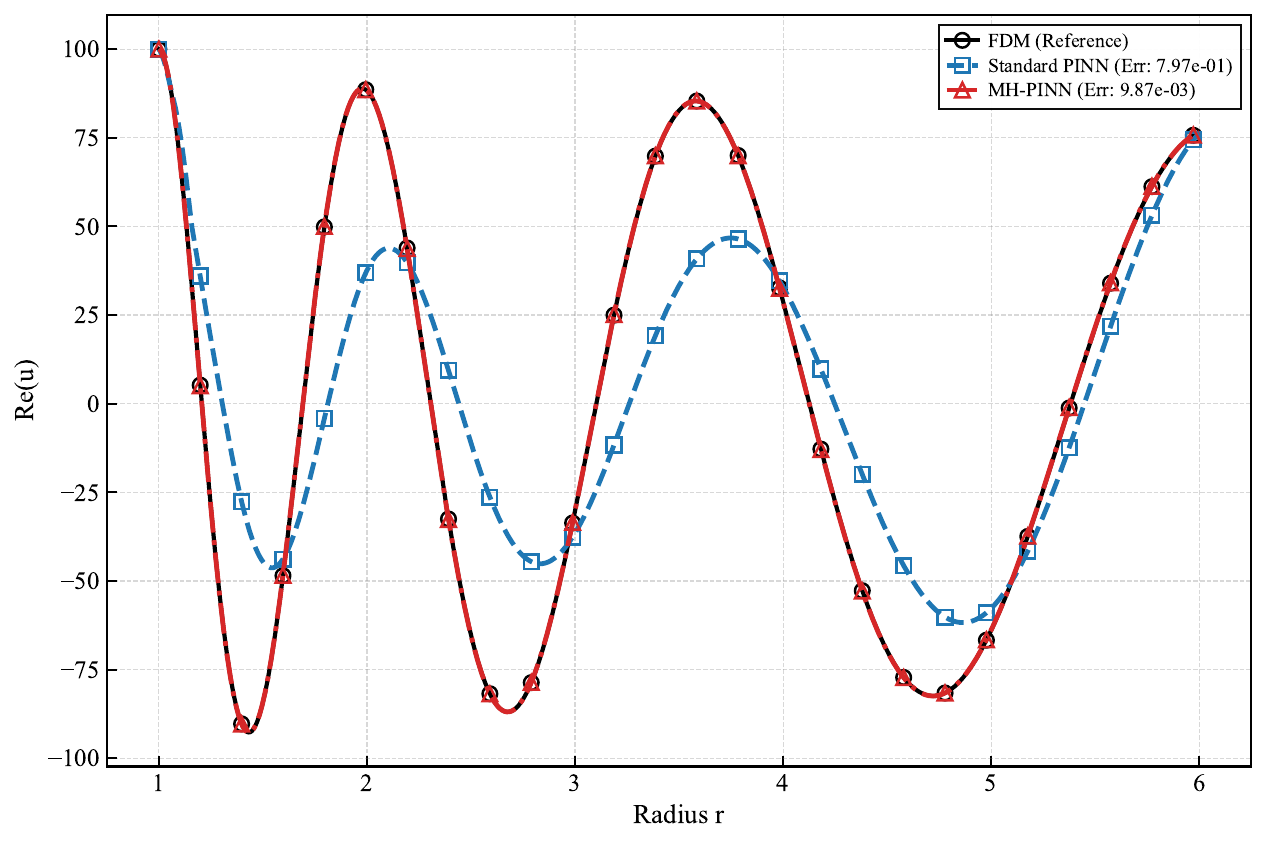}
    \caption{Comparison of calculation results of MH-PINN and FDM in the test area}
    \label{fig: varying coefficient Comparison}
\end{figure}

% 请确保导言区已包含 \usepackage{subcaption}

Figure \ref{fig: varying coefficient Comparison} shows a comparison between the predicted wavefield and the reference solution. The proposed method successfully captures the frequency modulation of the wave. A radial comparison of the sound pressure field clearly distinguishes the dense wave packets formed by rapid oscillations with high wavenumbers in the near field and the elongated wavelengths in the far field. The predicted solution agrees well with the FDM reference solution, accurately reproducing amplitude attenuation and phase shift caused by the variable medium. It is worth noting that without an integral phase term $\Phi(r)$ in the assumptions, the standard PINN algorithm is often affected by ``spectral bias,'' leading to non-convergence. Our method effectively mitigates spectral bias and significantly simplifies the optimization difficulty of the problem.

\subsection{Acoustic Scattering Problem}

This section examines the acoustic scattering problem in unbounded media. Unlike radiation problems, scattering problems involve strong coupling between incident and scattered waves, exhibiting distinctly different physical interference characteristics at different frequencies. We first verify the broadband adaptability of our method using the classic cylindrical scattering benchmark problem, then extend it to non-circular scatterers with complex geometries, and finally explore the potential of this framework in three-dimensional problems.

\subsubsection{Two-dimensional cylindrical benchmark test}

The physical model considers a soft acoustic cylinder of radius $R_{in} = 1$ in an infinite domain, subjected to a monochromatic incident plane wave $u_{inc}(\mathbf{x}) = e^{ikx}$ propagating along the positive $x$-axis. We aim to solve for the scattered field $u_s$. 

For a sound-soft circular cylinder, the exact analytical solution for the scattered acoustic pressure field can be expressed using the infinite Mie series expansion:
\begin{equation}
    u_{exact}(r, \theta) = - \sum_{n=0}^{\infty} \epsilon_n i^n \frac{J_n(kR_{in})}{H_n^{(1)}(kR_{in})} H_n^{(1)}(kr) \cos(n\theta)
    \label{eq:exact_scattering}
\end{equation}
where $\epsilon_n$ is the Neumann factor ($\epsilon_0=1$, and $\epsilon_n=2$ for $n > 0$), $J_n$ is the Bessel function of the first kind of order $n$, and $H_n^{(1)}$ is the Hankel function of the first kind of order $n$. This analytical solution serves as the ground truth for evaluating the computational accuracy of the proposed framework.

To evaluate the MH-PINN framework, we compare it against a baseline PINN that utilizes a traditional truncation strategy.

\textbf{Model Setup and Optimization:} The baseline PINN restricts the computational domain by introducing a first-order absorbing boundary condition (ABC) at an artificial truncation boundary $\Gamma_{out}$ ($r = 6$). It optimizes a soft-constraint loss function comprising weighted PDE residuals, inner boundary conditions, and ABC residuals.

In contrast, MH-PINN significantly simplifies the parameter requirements. Across the entire wideband test ($k \in [1, 10]$), MH-PINN maintains a unified, lightweight fully connected neural network architecture (4 hidden layers, 64 neurons per layer) with only $N_{col} = 12,000$ Latin hypercube sampling points generated within the coordinate-mapped domain $\xi \in [-1, 1), \theta \in [0, 2\pi)$. Both models employ a two-stage optimization strategy (Adam followed by L-BFGS) and are trained on an NVIDIA RTX 4060 GPU using PyTorch.

\begin{figure}[htbp] % 如果你之前用了 float 宏包，这里也可以换成 [H]
    \centering

    % 左侧子图 (a)：误差分析
    \begin{subfigure}[b]{0.48\textwidth}
        \centering
        \includegraphics[width=\textwidth]{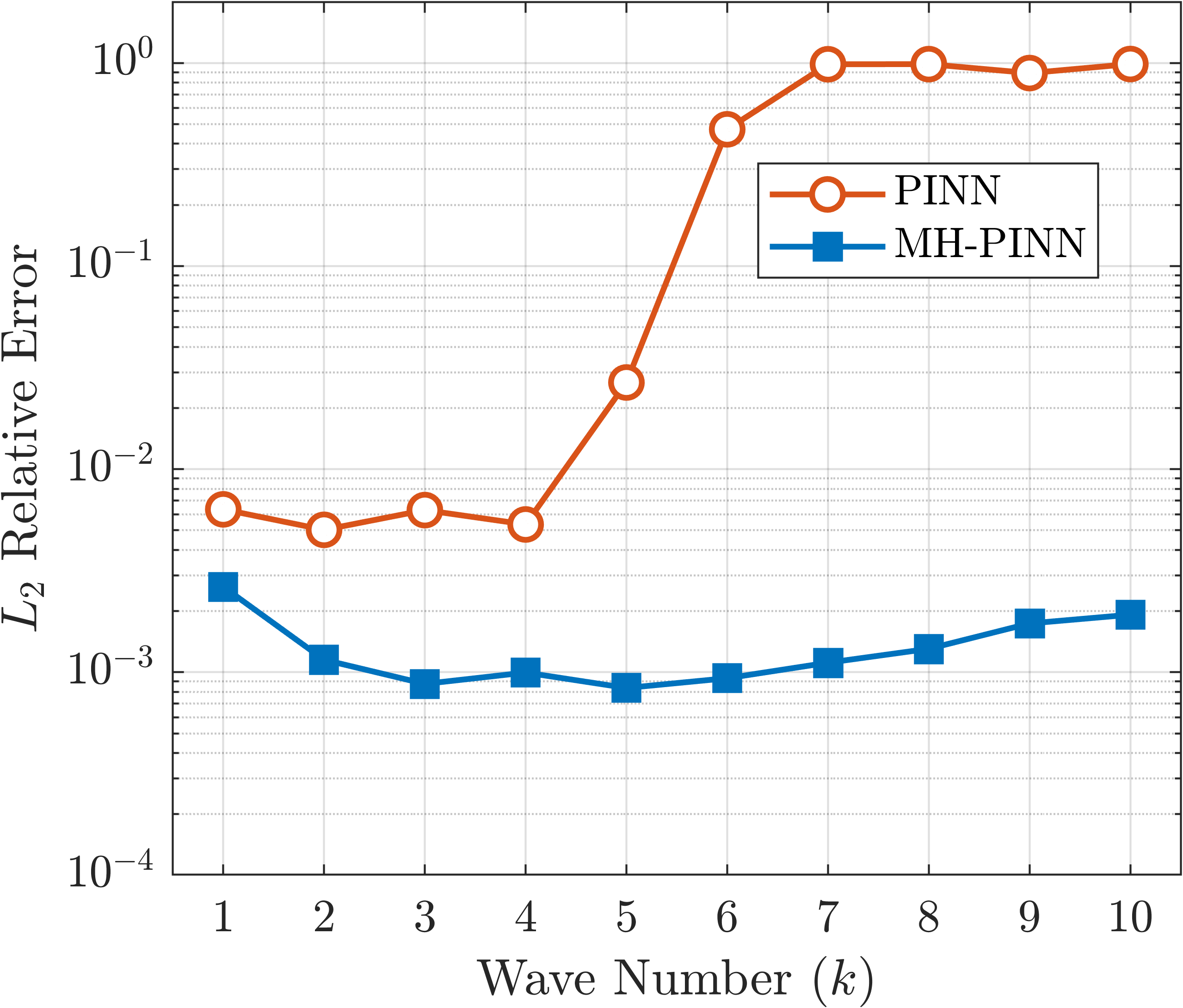}
        \caption{Computational Error Analysis}
        \label{fig:scatter_error}
    \end{subfigure}
    \hfill % hfill 会自动在两个子图之间插入适当的空白间距
    % 右侧子图 (b)：效率分析
    \begin{subfigure}[b]{0.48\textwidth}
        \centering
        \includegraphics[width=\textwidth]{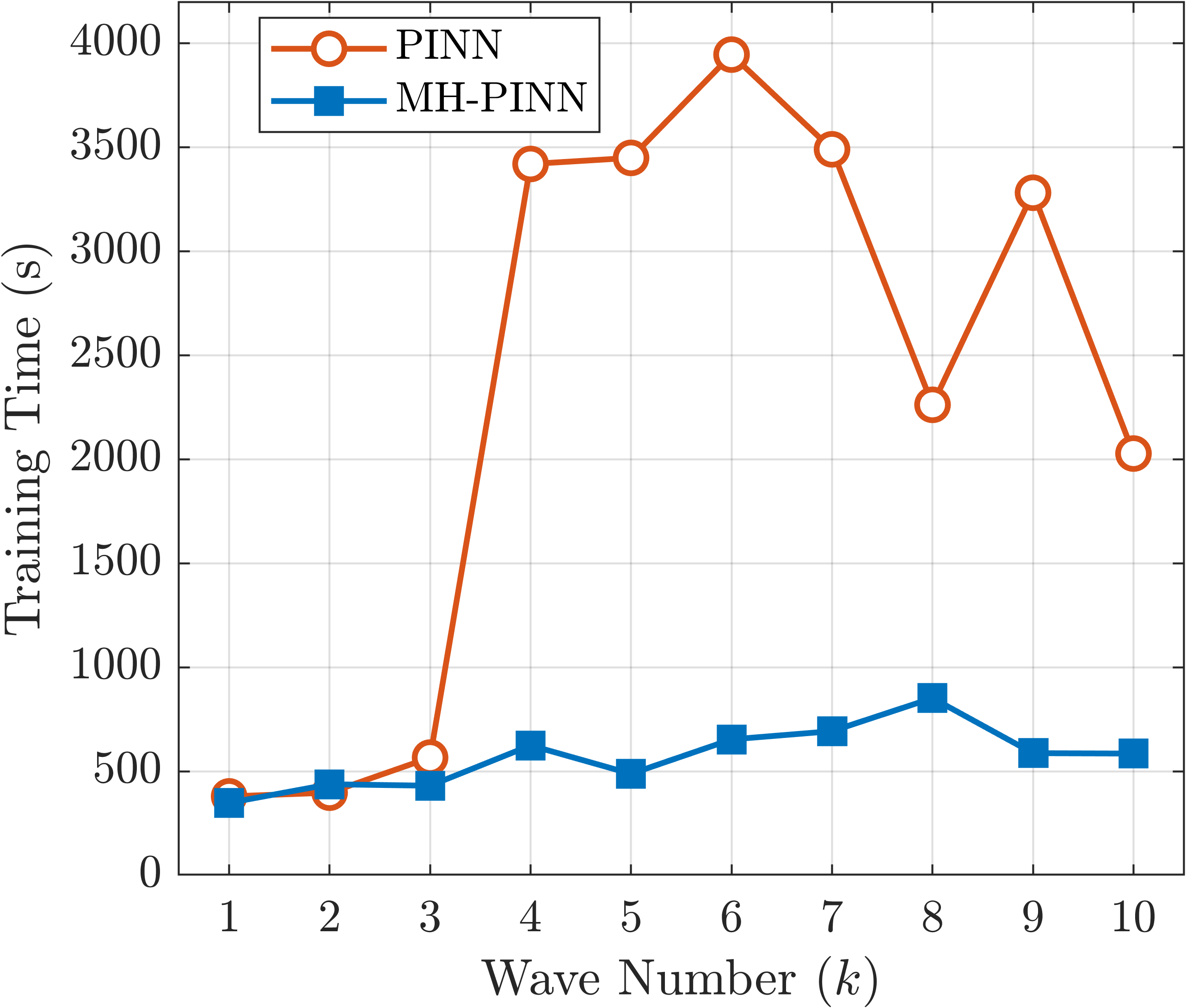}
        \caption{Computational Efficiency Analysis}
        \label{fig:scatter_efficiency}
    \end{subfigure}
    
% 整个大图的图注
    \caption{Comparison of computational performance between standard PINN and MH-PINN across different wavenumbers: (a) Error analysis; (b) Efficiency analysis.}
    \label{fig:performance_comparison}
\end{figure}

\textbf{Results and Discussion:} Figure \ref{fig:performance_comparison} compares the computational accuracy and efficiency of the two methods. Statistical results show that MH-PINN exhibits excellent frequency robustness: across the entire tested frequency band, the relative error remains consistently low, ranging from $\mathcal{O}(10^{-4})$ to $\mathcal{O}(10^{-3})$. Even at high frequencies of $k = 10$, the error is controlled to the order of $\mathcal{O}(10^{-3})$. 

Typically, as the wave number increases, the drastic oscillations in the solution can cause the solution space of the loss function to become extremely non-convex, making it easy for the optimizer to get trapped in local minima. However, thanks to the physical hard constraint assumptions in MH-PINN that partially decouple the frequency dependence of the problem, its error curve does not grow exponentially with the wave number $k$, averting the divergence seen in the baseline caused by high-frequency non-convexity. Furthermore, the training time of MH-PINN is significantly lower and more stable, verifying the high efficiency of this constraint-free optimization framework.

\subsubsection{Two-dimensional geometric generalization test}

After verifying the accuracy of circular boundaries, we further tested MH-PINN's ability to handle more general geometries. Specifically, we examined elliptical cylinders with geometric anisotropy, as well as square and hexagonal scatterers. For these non-circular geometries where the boundary radius $r_b(\theta)$ varies with angle, we directly applied the inverse factor correction strategy proposed in Section 3.2.2 to construct geometrically adaptive boundary coefficients $A(\theta)$ to achieve precise hard constraints on the Dirichlet boundary conditions. To systematically evaluate the algorithm's broadband adaptability, we conducted numerical simulations under typical operating conditions with three different wavenumbers: $k=1$, $3$, and $5$.

\begin{figure}[H]
    \centering
    \begin{subfigure}[b]{0.32\textwidth}
        \centering
        \includegraphics[width=\textwidth]{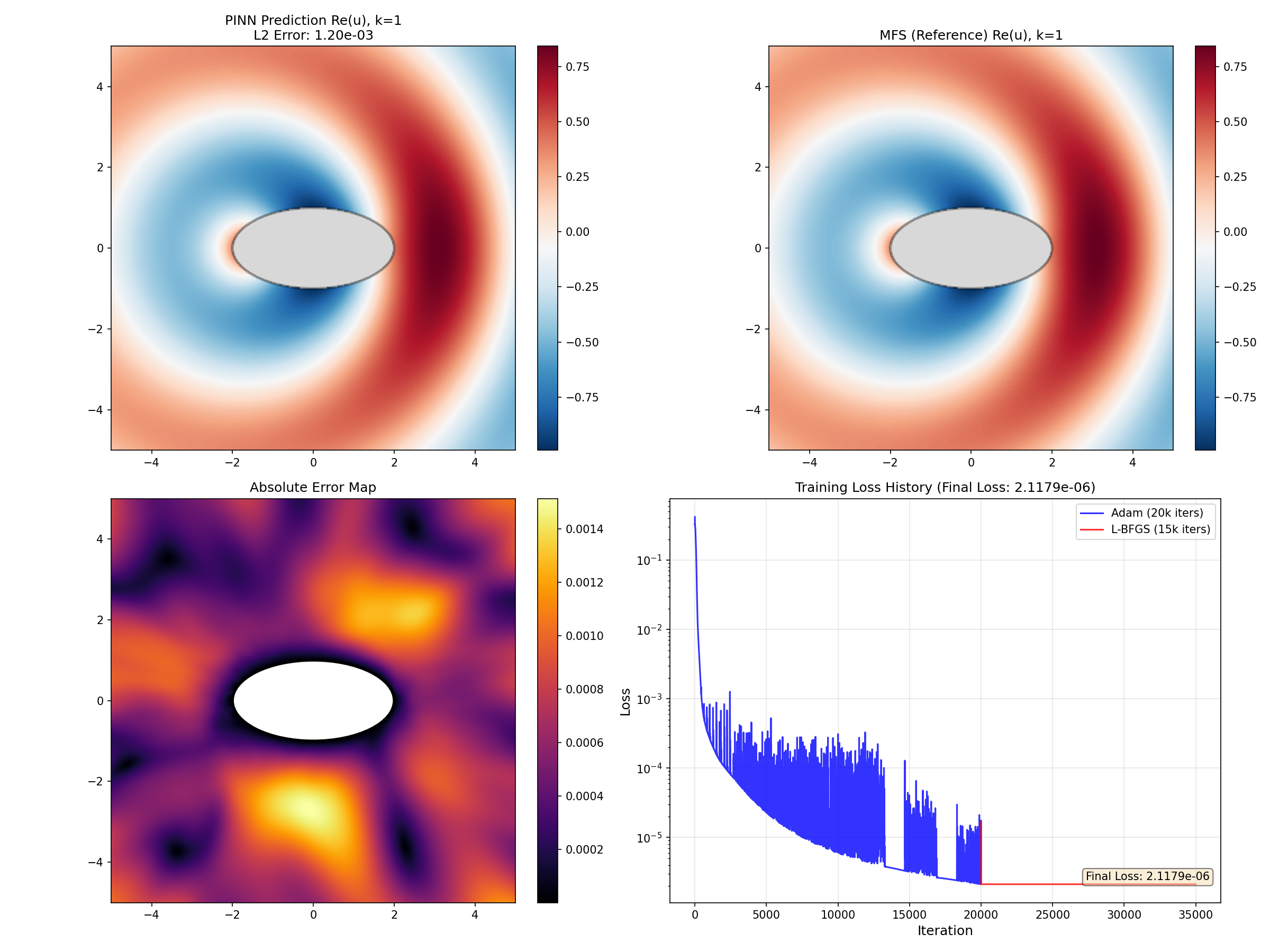} 
        \caption{$k=1$}
        \label{fig:elliptical_k1}
    \end{subfigure}
    \hfill % 水平弹性间距，让三张图均匀分布 
    \begin{subfigure}[b]{0.32\textwidth}
        \centering
        \includegraphics[width=\textwidth]{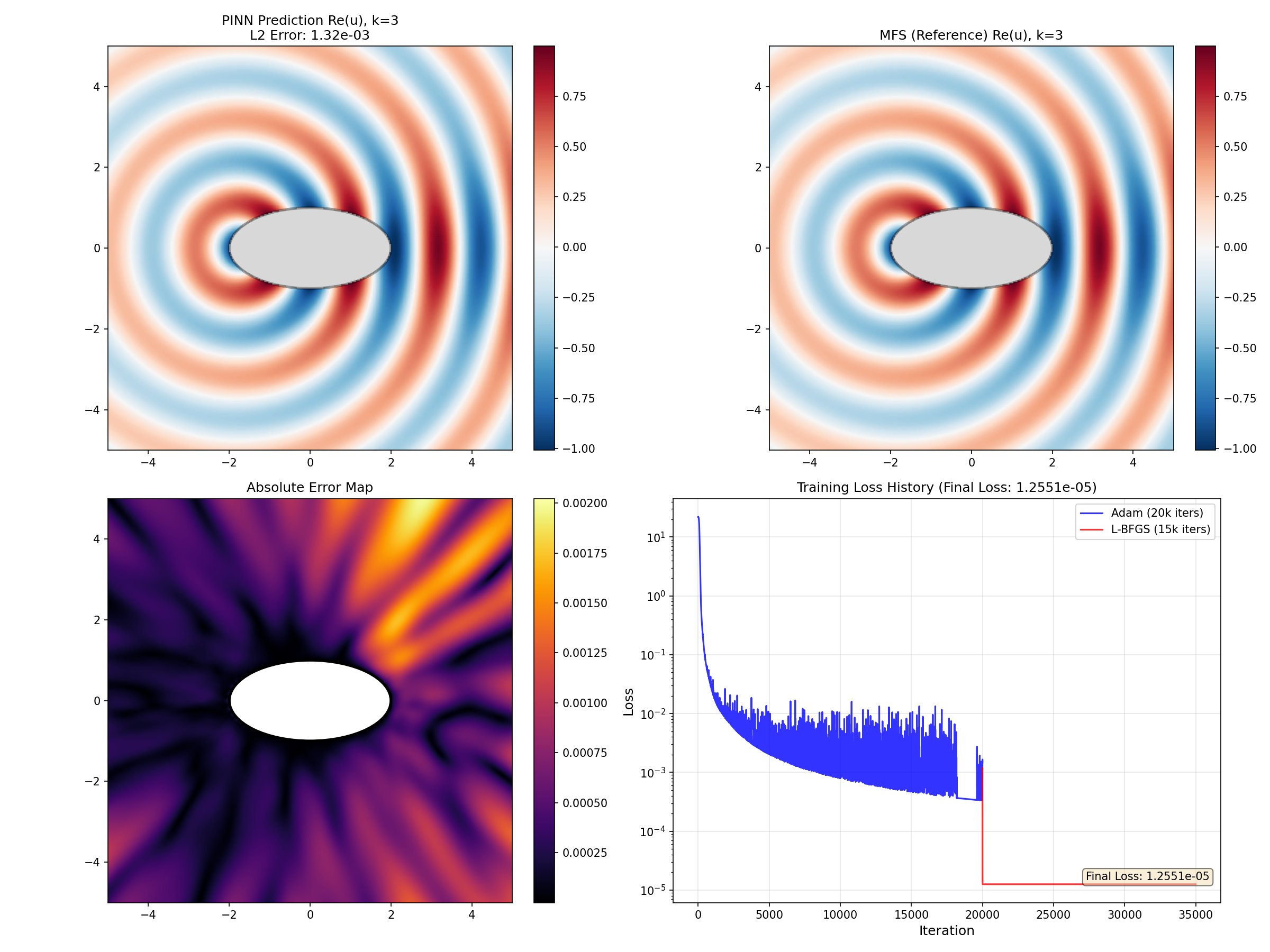} 
        \caption{$k=3$}
        \label{fig:elliptical_k3}
    \end{subfigure}
    \hfill % 水平弹性间距 
    \begin{subfigure}[b]{0.32\textwidth}
        \centering
        \includegraphics[width=\textwidth]{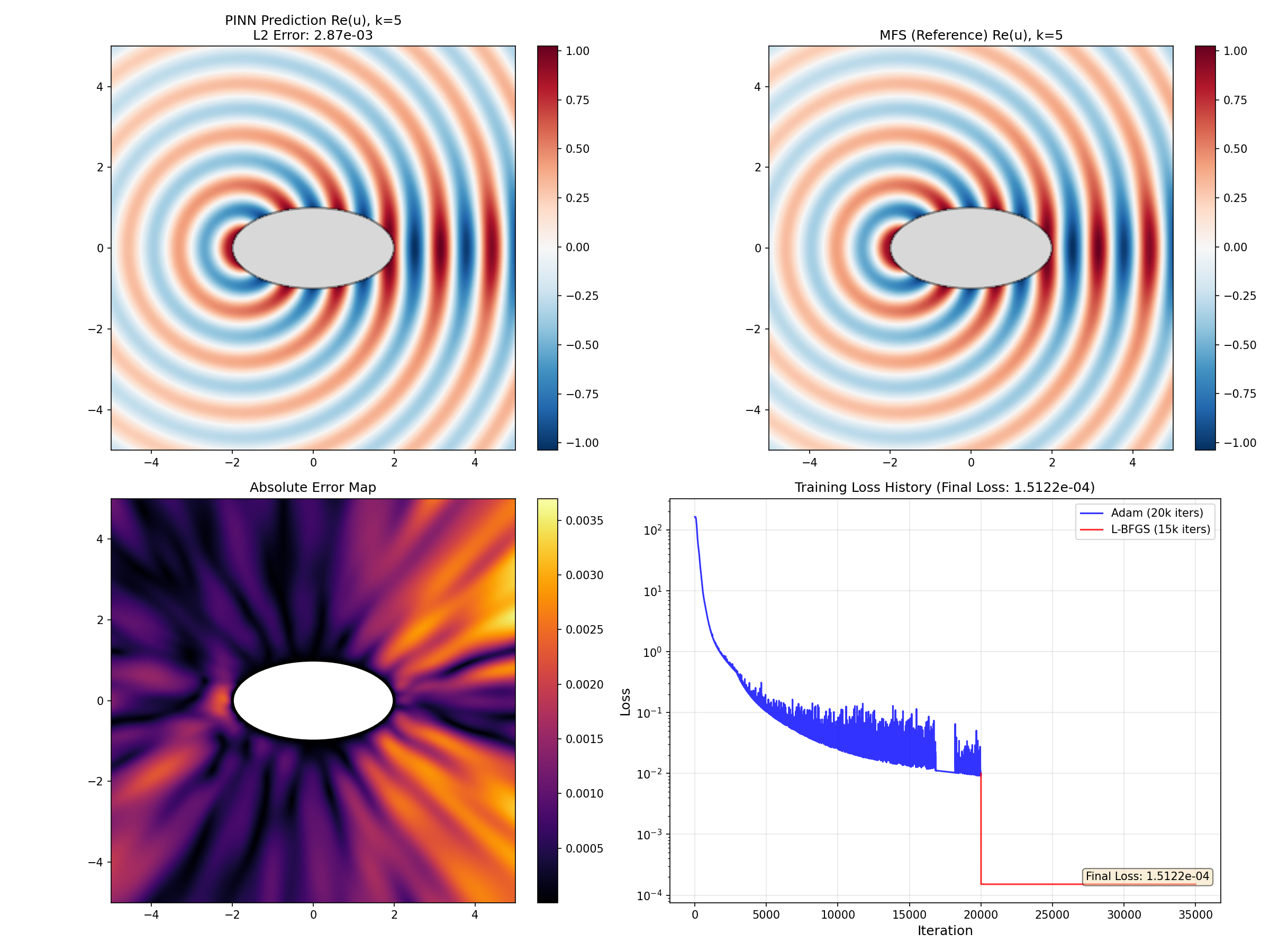} 
        \caption{$k=5$}
        \label{fig:elliptical_k5}
    \end{subfigure}
    \caption{Comparison of the real sound pressure part of MH-PINN and the MFS reference solution for different wavenumbers (elliptical boundary).}
    \label{fig:elliptical_scatterer}
\end{figure}

\begin{figure}[H]
    \centering
    \begin{subfigure}[b]{0.32\textwidth}
        \centering
        \includegraphics[width=\textwidth]{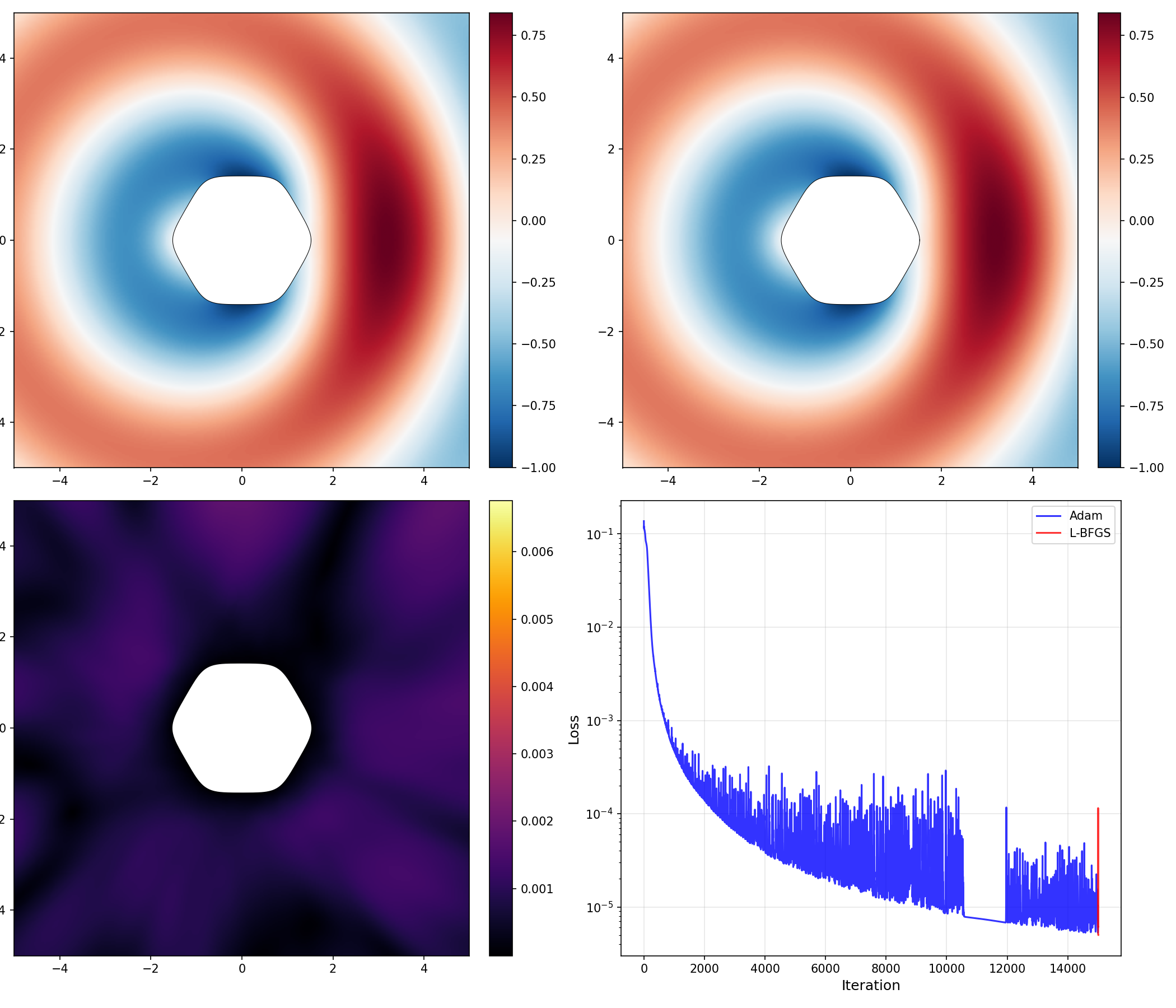}
        \caption{$k=1$}
        \label{fig:hexagonal_k1}
    \end{subfigure}
    \hfill
    \begin{subfigure}[b]{0.32\textwidth}
        \centering
        \includegraphics[width=\textwidth]{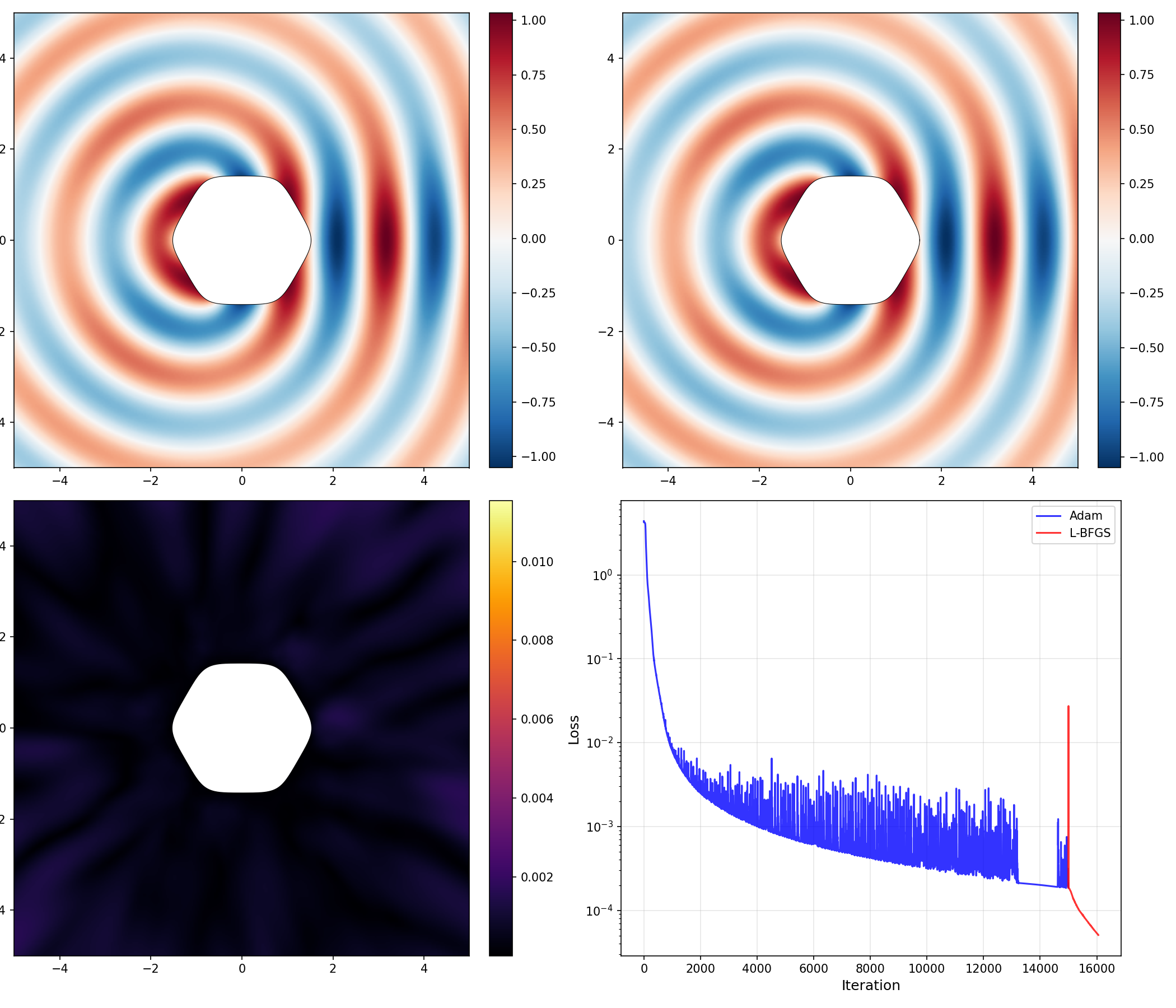}
        \caption{$k=3$}
        \label{fig:hexagonal_k3}
    \end{subfigure}
    \hfill
    \begin{subfigure}[b]{0.32\textwidth}
        \centering
        \includegraphics[width=\textwidth]{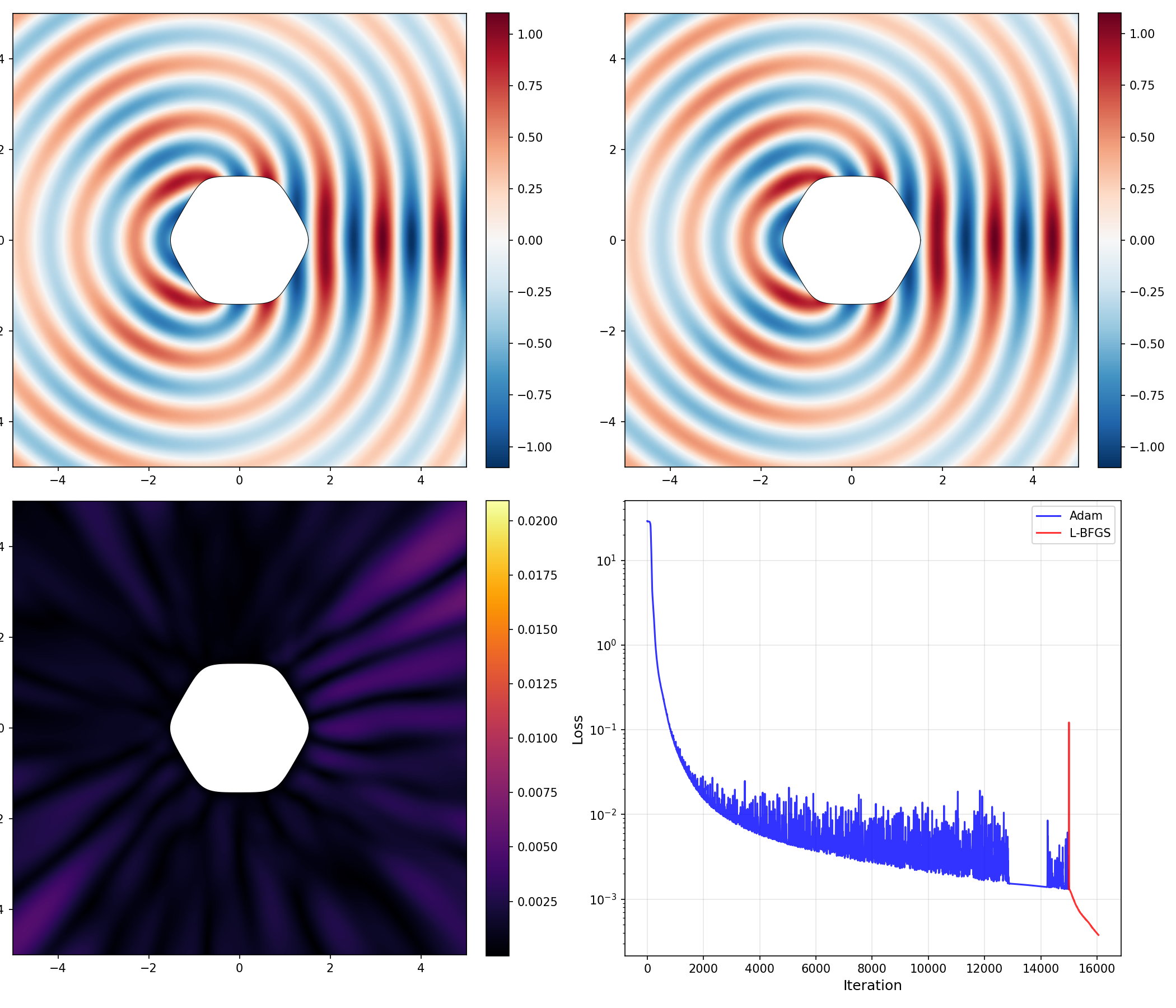}
        \caption{$k=5$}
        \label{fig:hexagonal_k5}
    \end{subfigure}
    \caption{Comparison of the real sound pressure part of MH-PINN and the MFS reference solution for different wavenumbers (hexagonal boundary).}
    \label{fig:hexagonal_scatterer}
\end{figure}

\begin{figure}[H]
    \centering
    
    % --- k=1 组 ---
    \begin{minipage}{\textwidth}
        \centering
% 18a1 (云图): 设置稍宽，0.65\textwidth
        \begin{subfigure}[b]{0.65\textwidth}
            \centering
            \includegraphics[width=\textwidth]{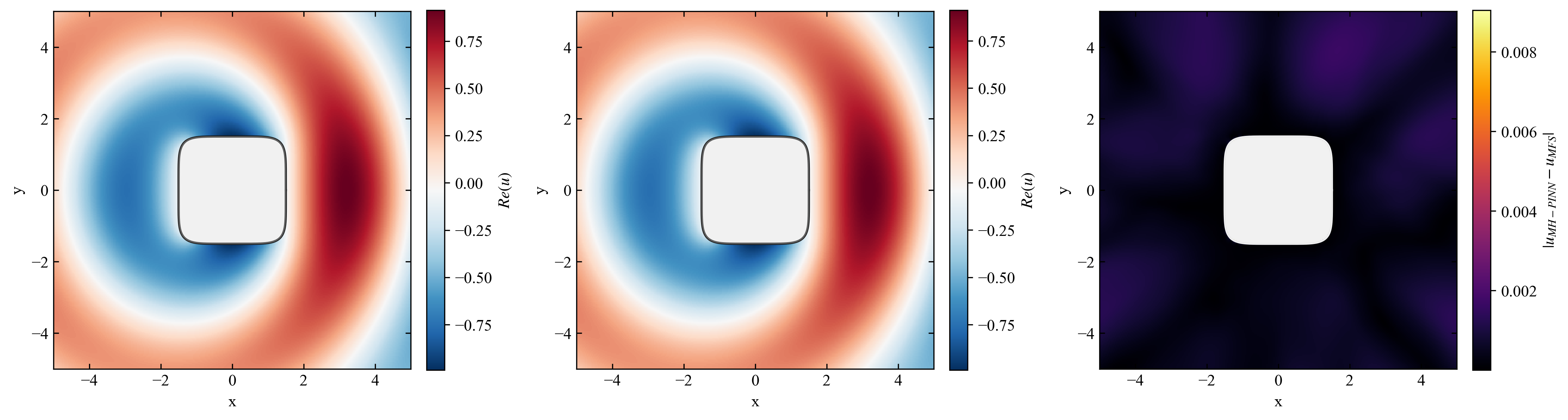}
            \caption{$k=1$: Sound pressure }
            \label{fig:square_k1_field}
        \end{subfigure}
        \hfill % 推到最右侧
        % 18a2 (Loss图): 设置稍窄，0.30\textwidth 
        \begin{subfigure}[b]{0.30\textwidth}
            \centering
            \includegraphics[width=\textwidth]{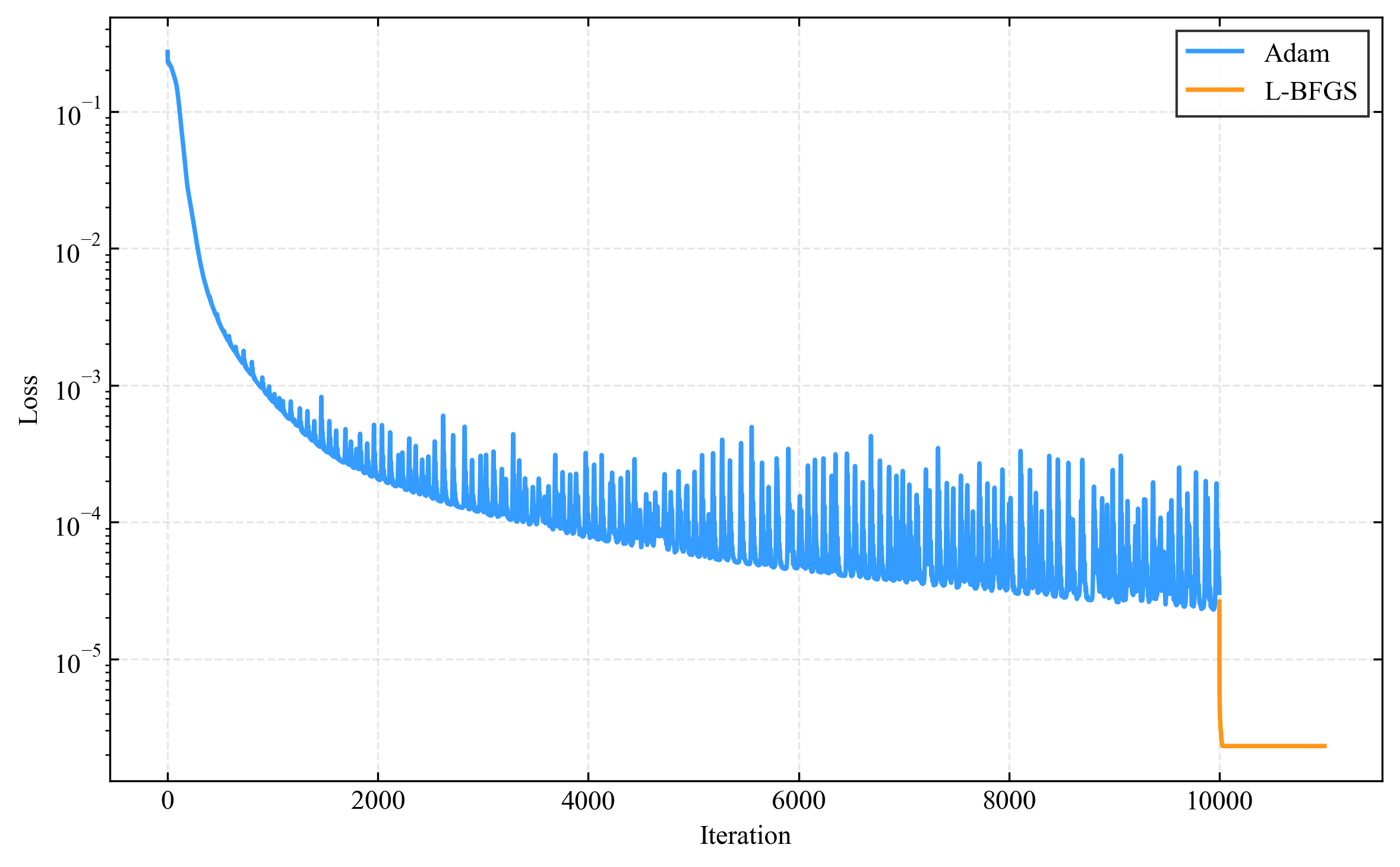}
            \caption{$k=1$: Loss evolution}
            \label{fig:square_k1_loss}
        \end{subfigure}
    \end{minipage}
    
    \vspace{0.4cm} % 组间距
 
    % --- k=3 组 ---
    \begin{minipage}{\textwidth}
        \centering
% 18b1 (云图)
        \begin{subfigure}[b]{0.65\textwidth}
            \centering
            \includegraphics[width=\textwidth]{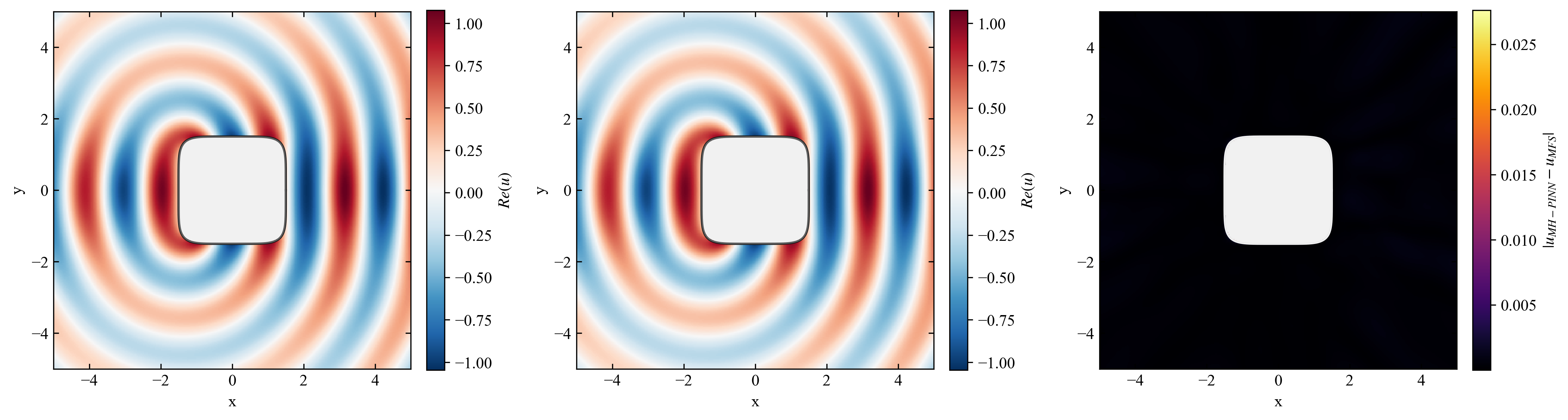}
            \caption{$k=3$: Sound pressure}
            \label{fig:square_k3_field}
        \end{subfigure}
        \hfill
% 18b2 (Loss图)
        \begin{subfigure}[b]{0.30\textwidth}
            \centering
            \includegraphics[width=\textwidth]{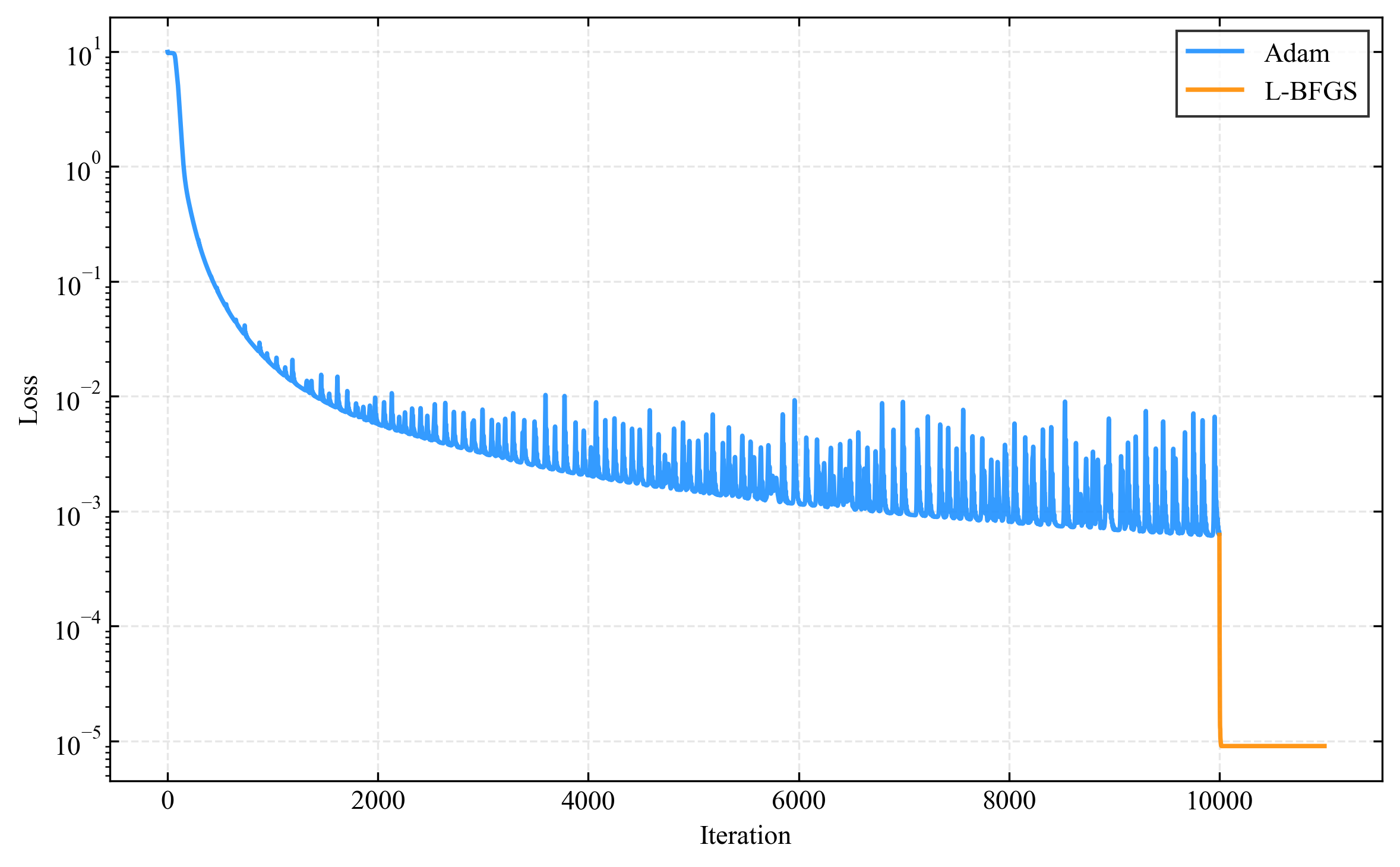}
            \caption{$k=3$: Loss evolution}
            \label{fig:square_k3_loss}
        \end{subfigure}
    \end{minipage}
    
    \vspace{0.4cm} % 组间距
 
    % --- k=5 组 ---
    \begin{minipage}{\textwidth}
        \centering
% 18c1 (云图)
        \begin{subfigure}[b]{0.65\textwidth}
            \centering
            \includegraphics[width=\textwidth]{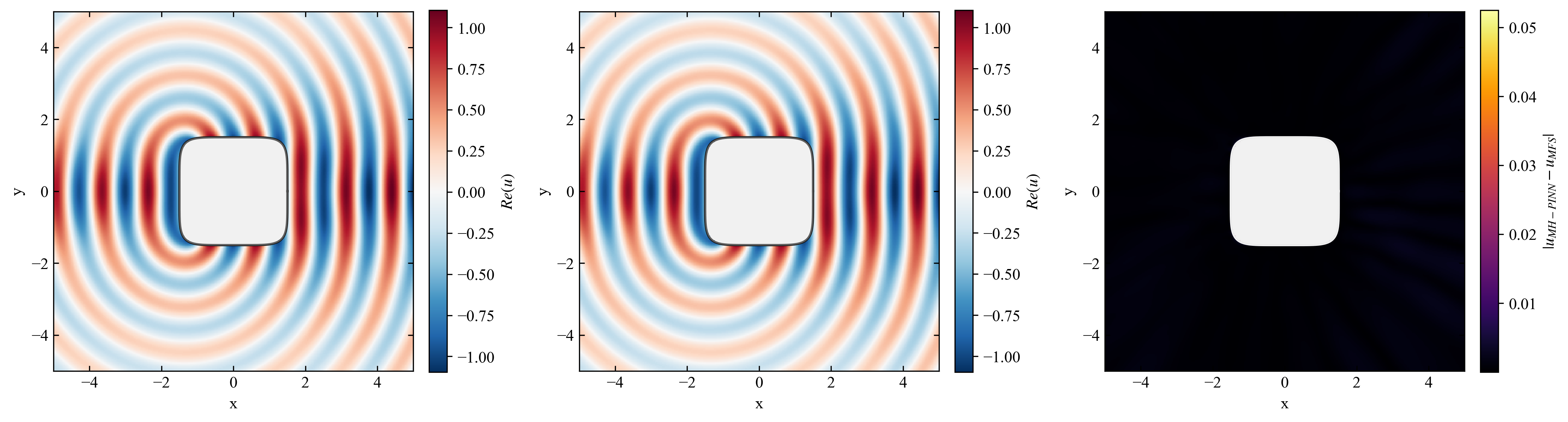}
            \caption{$k=5$: Sound pressure}
            \label{fig:square_k5_field}
        \end{subfigure}
        \hfill
% 18c2 (Loss图)
        \begin{subfigure}[b]{0.30\textwidth}
            \centering
            \includegraphics[width=\textwidth]{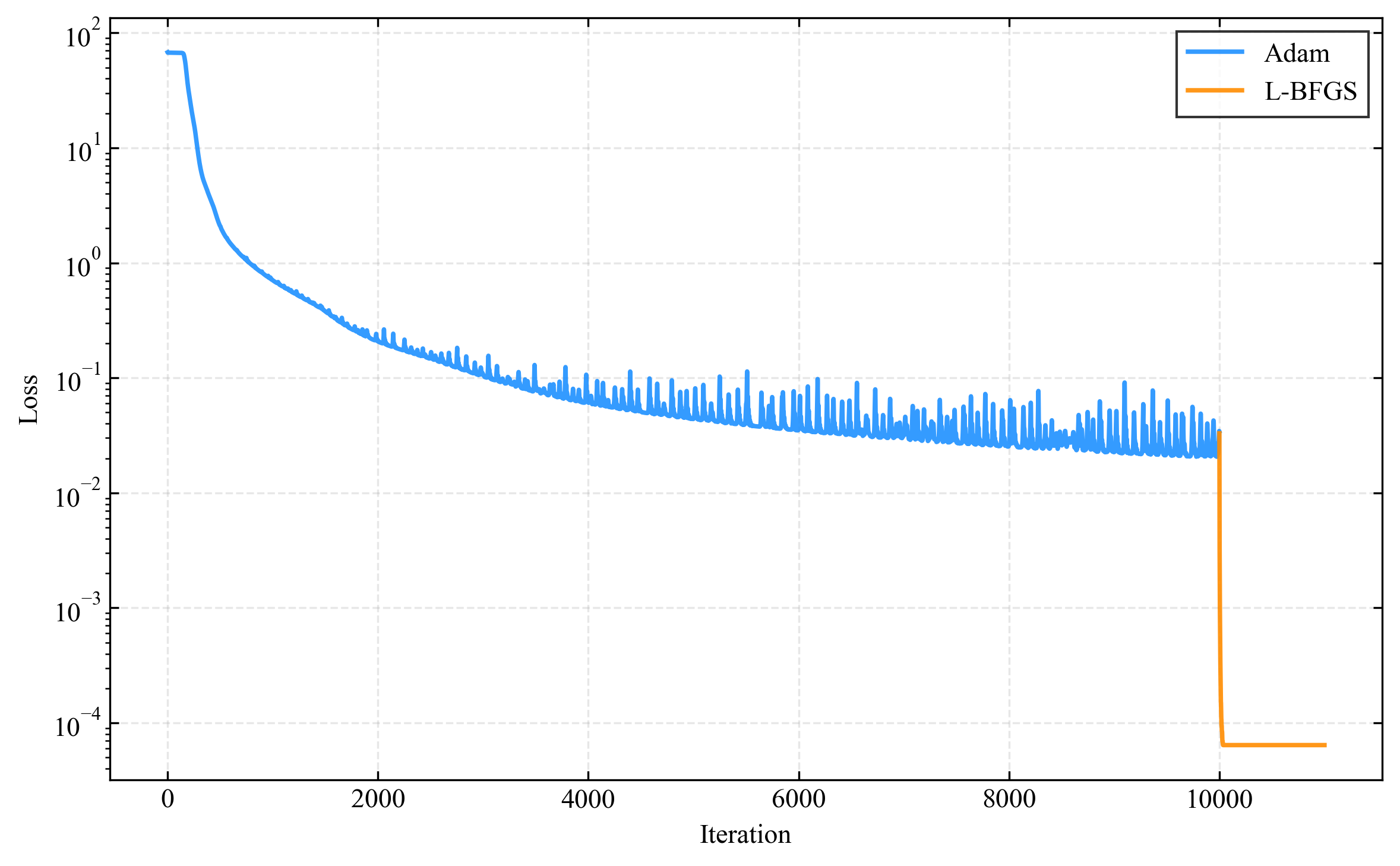}
            \caption{$k=5$: Loss evolution}
            \label{fig:square_k5_loss}
        \end{subfigure}
    \end{minipage}
    
    \caption{Numerical results for the square boundary scatterer at different wavenumbers. For each case, the left panel displays the field distribution (MH-PINN solution, MFS reference, and absolute error), and the right panel displays the training loss history. The results demonstrate the ability of the proposed MH-PINN to accurately handle scatterers with sharp corners at higher frequencies.}
    \label{fig:square_scatterer}
\end{figure}

Figures \ref{fig:elliptical_scatterer}, \ref{fig:hexagonal_scatterer}, and \ref{fig:square_scatterer} show the real part contour plots of sound pressure and the distribution of absolute error for elliptical, hexagonal, and square scatterers at different wavenumbers, respectively. The reference solution is calculated using the Method of Fundamental Solutions.

For the elliptical scatterer, despite its geometry breaking rotational symmetry, MH-PINN accurately reconstructs the scattered wave field. The predicted solution is visually highly consistent with the MFS reference solution, indicating that coordinate mapping and inverse factor correction successfully solve the mesh adaptation problem for anisotropic geometry. The absolute error contour plots show that the error is mainly concentrated in the minimal neighborhood of the corners,and the overall $L_2$ relative error of the test cases is on the order of magnitude or below $\mathcal{O}(10^{-3})$.

In summary, MH-PINN demonstrates excellent numerical accuracy and convergence stability when handling anisotropic ellipses or polygons of other shapes. This indicates that the framework possesses good geometric generalization ability, providing a unified, efficient, and mesh-reconstruction-free computational scheme for unbounded acoustic simulation of scatterers of general shapes.

% ---------------- Figure 9: Elliptical boundary ----------------

% ---------------- Figure 10: Hexagonal boundary ----------------

% ---------------- Figure 11: Square boundary scatterer (Complete Grid) ----------------

\subsubsection{Three-dimensional sphere benchmark test}
This section applies the proposed MH-PINN method to the problem of acoustic wave scattering in a three-dimensional unbounded domain. To evaluate the computational accuracy of the proposed method, we first conduct a benchmark test using a standard acoustic soft sphere scattering example with an exact analytical solution.

We first consider the 3D acoustic scattering by a sound-soft sphere of radius $r_0 = 1.0$ centered at the origin. Assuming an incident plane wave propagating along the positive $x$-axis, denoted as $u^{inc} = e^{ikx}$, a homogeneous Dirichlet boundary condition is imposed on the surface of the scatterer, indicating a zero total acoustic pressure ($u^{scat} = -u^{inc}$). For such a standard spherical geometry, the exact analytical solution of the scattered field in the unbounded exterior domain can be rigorously expressed by the infinite Mie series expansion:
\begin{equation}
    u_{exact}(r, \theta) = -\sum_{n=0}^{\infty} i^n (2n+1) \frac{j_n(kr_0)}{h_n^{(1)}(kr_0)} h_n^{(1)}(kr) P_n(\cos\theta)
    \label{eq:exact_scattering_3d}
\end{equation}
where $j_n$ denotes the spherical Bessel function of the first kind of order $n$, $h_n^{(1)}$ is the spherical Hankel function of the first kind of order $n$, $P_n$ represents the Legendre polynomial of degree $n$, and $\theta$ is the polar angle measured from the positive $x$-axis (the direction of wave propagation). This exact closed-form solution serves as the rigorous ground truth for evaluating our 3D predictions.

Within the MH-PINN framework, the infinite exterior physical domain is mapped onto a finite computational domain $\xi \in [-1, 1)$ using a conformal coordinate transformation with a scaling parameter $L = 2.0$. To evaluate the model's performance in a highly oscillatory field, the dimensionless wavenumber is set to $k = 5.0$. Thanks to the introduction of the hard-constraint ansatz, the network output is forced to strictly satisfy both the physical boundary condition and the Sommerfeld radiation condition at infinity. Consequently, during the training process, no data points need to be sampled on the boundary or at the far-field truncation. We only randomly generate $N_{col} = 10,000$ collocation points within the interior of the computational domain to minimize the PDE residual of the Helmholtz equation.

\begin{figure}[H]
    \centering
    % 子图 (a) 实部
    \includegraphics[width=0.95\textwidth]{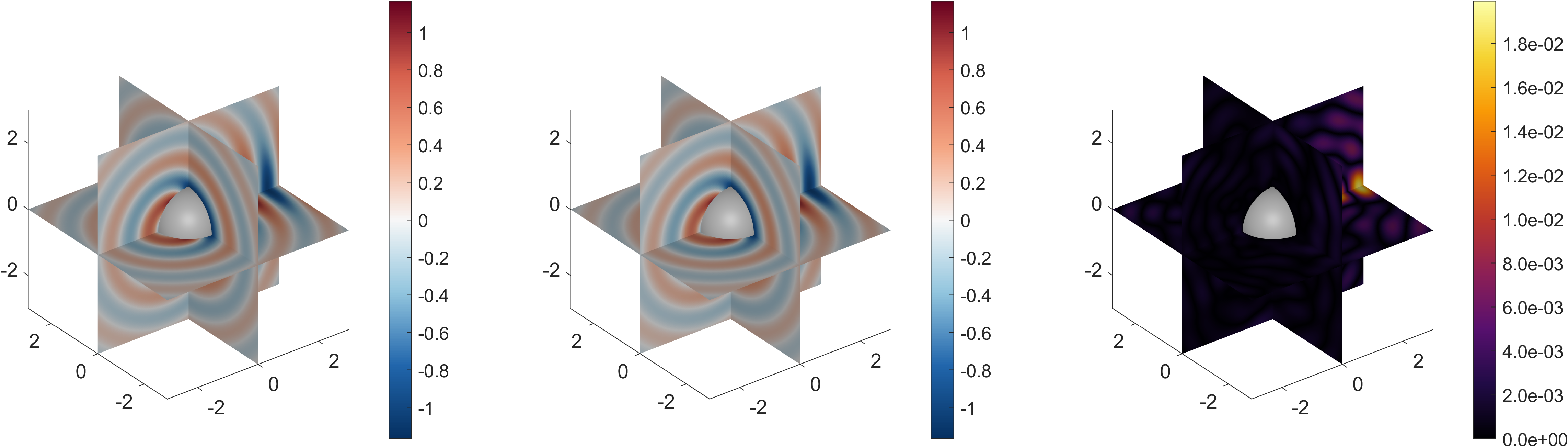} \\
    \vspace{0.1cm}
    \small{(a) Real part} \\
    \vspace{0.4cm}
    % 子图 (b) 虚部
    \includegraphics[width=0.95\textwidth]{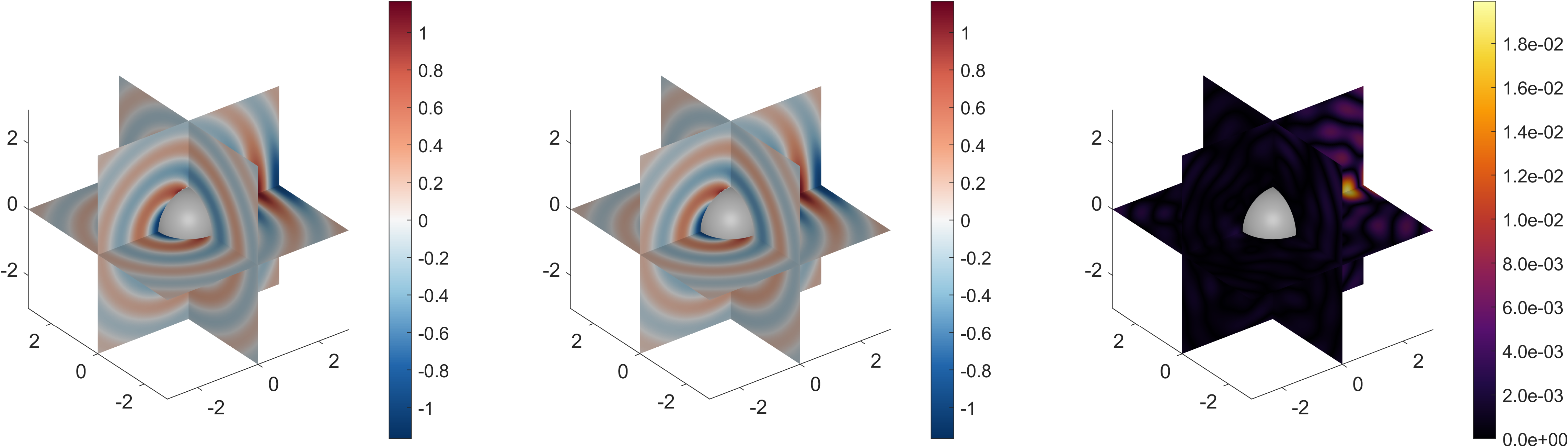} \\
    \vspace{0.1cm}
    \small{(b) Imaginary part}
    
    \caption{Acoustic scattering by a sound-soft sphere ($k=5.0$): spatial distribution of the scattered acoustic pressure field on the orthogonal central cross-sections ($x=0$, $y=0$, and $z=0$). (a) Real part; (b) Imaginary part. In both (a) and (b), the left, middle, and right panels display the MH-PINN prediction, the exact analytical solution, and the absolute error, respectively.}
    \label{fig:sphere_scattering}
\end{figure}

Figure \ref{fig:sphere_scattering} illustrates the spatial distribution characteristics of the scattered acoustic pressure field on three orthogonal central cross-sections ($x=0$, $y=0$, and $z=0$) after network training. The figure provides a comparison among the MH-PINN predictions, the exact analytical solutions, and their absolute errors, comprehensively covering both the real and imaginary parts of the complex acoustic field. It can be intuitively observed that the MH-PINN predictions are in excellent agreement with the exact analytical solutions. Furthermore, the absolute error maps indicate that the pointwise errors across the entire unbounded computational domain are maintained at an extremely low magnitude. It is particularly noteworthy that the absolute error strictly converges to zero at the spherical boundary surface ($r \to r_0$). This phenomenon strongly demonstrates that the proposed analytical hard-constraint strategy not only completely eliminates the difficulty of weight tuning caused by soft constraints in conventional PINNs but also fundamentally guarantees the ultimate approximation accuracy of the numerical solution near the physical boundary.

% --- Figure for Sphere ---

\subsubsection{Three-dimensional geometric generalization test}

After thoroughly verifying the basic accuracy, we further extend MH-PINN to a three-dimensional acoustic soft ellipsoid scattering example. Since obtaining closed-form analytical solutions in the ellipsoidal coordinate system is extremely difficult and subject to truncation errors, this section uses the high-fidelity results calculated by the Method of Fundamental Solutions as the reference truth for comparison and evaluation.

We set the lengths of the three semi-axis of the ellipsoid to be $a=2.0$, $b=1.0$, and $c=1.0$, and the external environment is also excited by plane waves, with wave number $k=5.0$. In this example, the mapping function is naturally generalized to a nonlinear coordinate transformation based on the local polar radius of the ellipsoid:
\begin{equation}
    r(\theta, \phi) = r_b(\theta, \phi) + L \frac{1+\xi}{1-\xi}
\end{equation}

This allows the physical boundary to be flattened and mapped to a regular computational plane of $\xi = -1$, thus avoiding the mesh generation problem under complex geometry.

% 调整位置：将 Figure 13 紧跟在第一次引用的段落下方，并适当缩小尺寸
\begin{figure}[htbp]
    \centering
    % 子图 (a)：实部，缩小宽度以节省垂直空间
    \begin{subfigure}{\textwidth}
        \centering
        \includegraphics[width=0.8\textwidth]{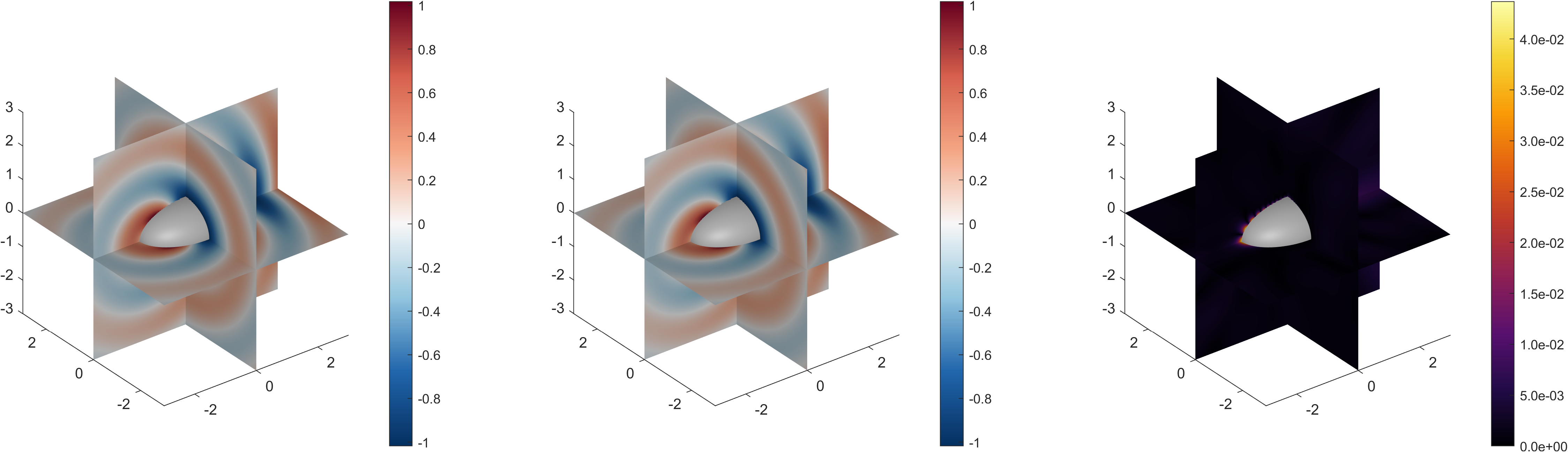}
        \caption{Real part}
        \label{fig:ellipsoid_real}
    \end{subfigure}
    
    \vspace{0.2cm} % 压缩子图之间的间距
    
    % 子图 (b)：虚部
    \begin{subfigure}{\textwidth}
        \centering
        \includegraphics[width=0.8\textwidth]{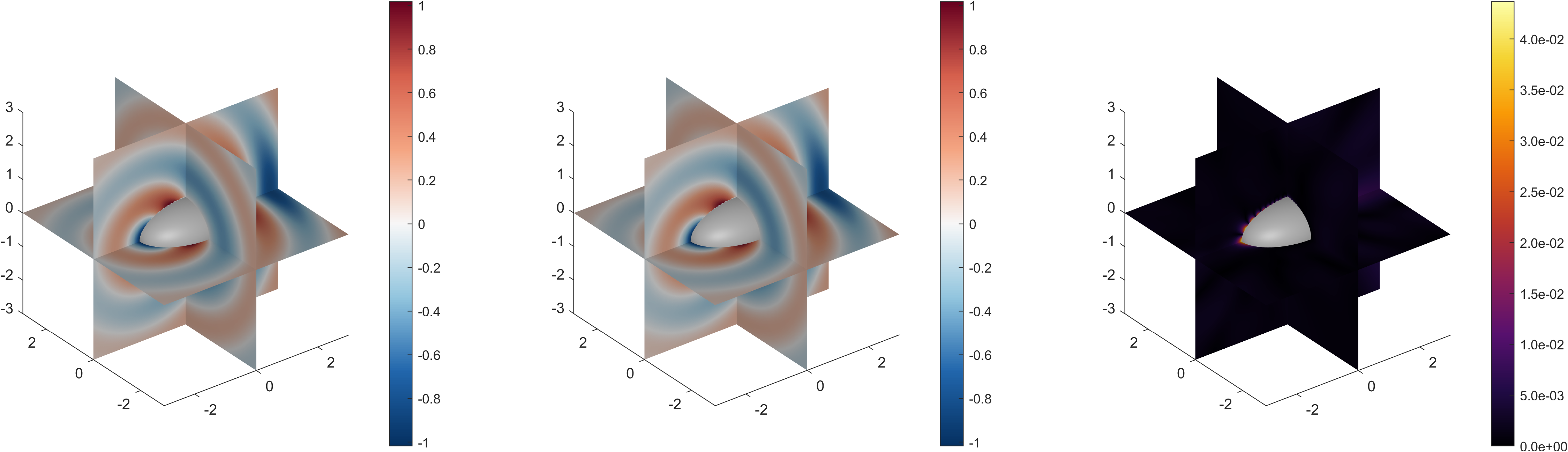}
        \caption{Imaginary part}
        \label{fig:ellipsoid_imag}
    \end{subfigure}
    
    \caption{Comparison results of the ellipsoidal scattering field on a three-dimensional slice.}
    \label{fig:ellipsoid_scattering}
\end{figure}

Figure \ref{fig:ellipsoid_scattering} shows the comparison results of the ellipsoidal scattering field on a three-dimensional slice. Despite the change in boundary shape, MH-PINN continues to demonstrate excellent predictive capability, with the predicted real and imaginary acoustic pressure fields showing high agreement with the reference solution obtained by the Method of Fundamental Solutions. The relative error of the real part L2 within the shown range is $6.04\times10^{-3}$. This example fully demonstrates that the accuracy of MH-PINN does not depend on regular geometry, and that the proposed adaptive coordinate mapping technique and hard-constrained network architecture can be seamlessly and efficiently generalized to various complex engineering three-dimensional acoustic scattering problems.

\subsection{Elastodynamic Scattering: SH Waves over a 2D Canyon}
\label{subsec:sh_wave_canyon}

To demonstrate the generalizability of the proposed MH-PINN framework to Neumann boundary conditions, we extend our application from acoustics to elastodynamics. We consider the classic problem of out-of-plane horizontally polarized shear (SH) wave scattering by a two-dimensional semi-circular canyon. 

As illustrated in previous elastodynamic studies, the semi-circular canyon with radius $a=1.0$ is located at the origin of an elastic half-space ($y \le 0$). The surface of the half-space and the canyon are traction-free, which translates to a homogeneous Neumann boundary condition ($\partial u / \partial n = 0$). To evaluate the frequency dependence, the dimensionless frequency $\eta$ is defined as:
\begin{equation}
    \eta = \frac{\omega a}{\pi c_s} = \frac{k a}{\pi}
\end{equation}
where $c_s$ is the shear wave velocity and $k$ is the wavenumber. In this example, we investigate the case of $\eta = 1$ (i.e., $k = \pi$). 

To apply the MH-PINN framework, we utilize the method of images to transform the half-space scattering problem into an equivalent full-space problem. Assuming the incident SH wave is a time-harmonic plane wave with an incident angle $\theta_{inc}$, its expression is given by:
\begin{equation}
    u_{inc}(x,y) = e^{ik(x\cos\theta_{inc} + y\sin\theta_{inc})}
\end{equation}
Because the traction on the flat ground surface ($y=0$, excluding the canyon part) must be identically zero, we introduce a mirror reflected wave $u_{ref}$ to satisfy this condition automatically:
\begin{equation}
    u_{ref}(x,y) = e^{ik(x\cos\theta_{inc} - y\sin\theta_{inc})}
\end{equation}
The superposition of these two waves constitutes the total background field $u^{(0)} = u_{inc} + u_{ref}$. By extending the semi-circular canyon into a full circular cavity in the full space, the scattering problem is elegantly simplified. Consequently, the boundary condition on the equivalent circular cavity ($r=a$) transforms into an inhomogeneous Neumann condition for the scattered field $u_{scat}$:
\begin{equation}
    \frac{\partial u_{scat}}{\partial r} \bigg|_{r=a} = -\frac{\partial u^{(0)}}{\partial r} \bigg|_{r=a} \triangleq g_N(\theta)
\end{equation}
Using the technique derived in Section \ref{subsec:hard_constraints}, this normal derivative constraint is strictly enforced as a hard architectural constraint. A mapping parameter $L=3.0$ and $20,000$ interior collocation points are utilized without any boundary penalty terms in the loss function.

\begin{figure}[htbp]
    \centering
    % 插入拼接后的大图。
    \includegraphics[height=0.85\textheight]{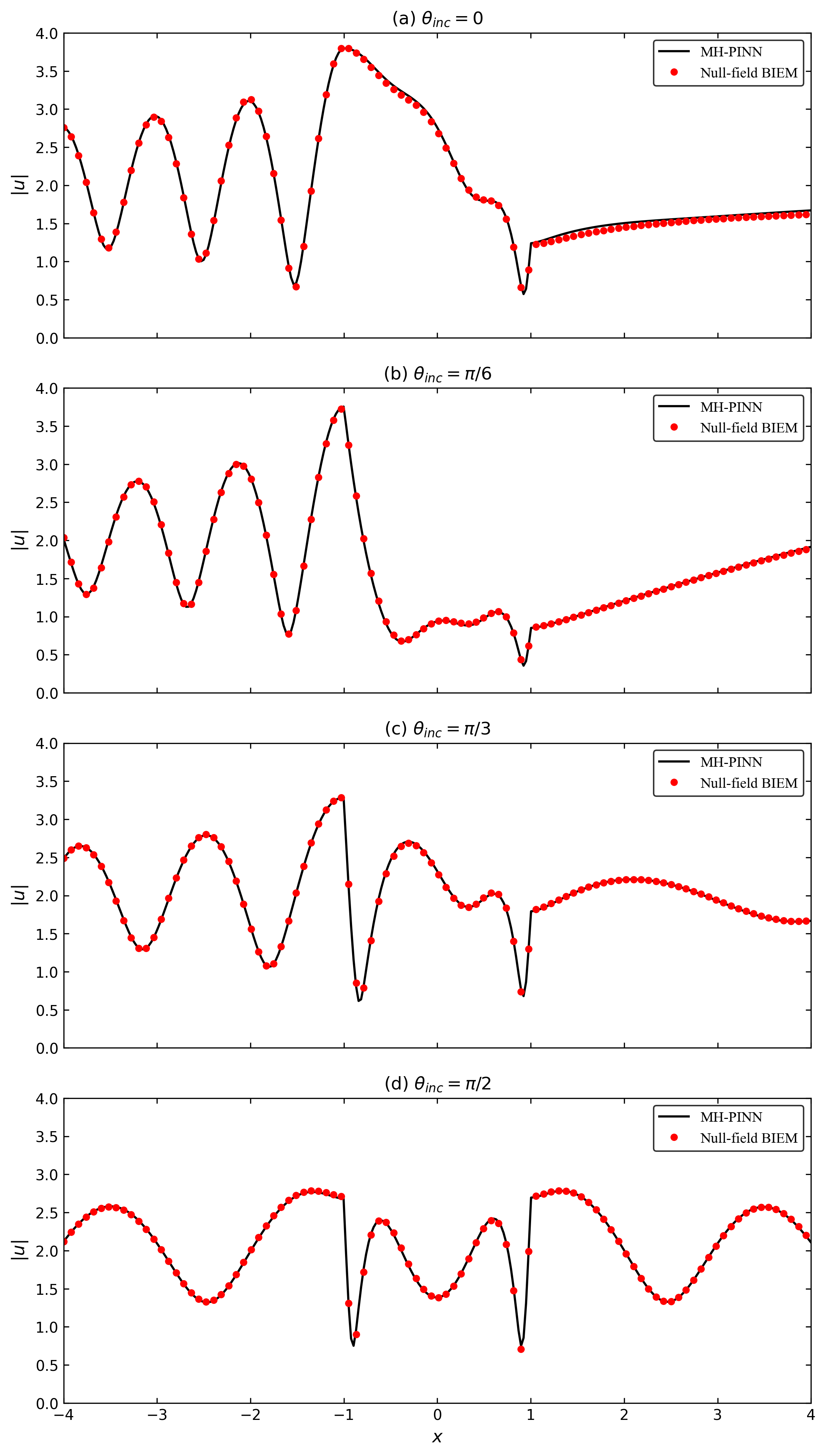} 
    
    % 稍微缩小图注与图片的距离，使其更紧凑
    \vspace{-0.2cm}
    \caption{Comparison of surface displacement amplitudes $|u|$ for SH wave scattering by a semi-circular canyon ($\eta=1$) at different incident angles. From top to bottom: (a) grazing incidence $\theta_{inc} = 0$, (b) oblique incidence $\theta_{inc} = \pi/6$, (c) oblique incidence $\theta_{inc} = \pi/3$, and (d) vertical incidence $\theta_{inc} = \pi/2$. The solid black lines represent the MH-PINN predictions using analytical compensation, while the red dots denote the Null-field BIEM reference solutions.}
    \label{fig:sh_wave_surface}
\end{figure}

We systematically investigate the influence of the incident angle $\theta_{inc}$ on the surface displacement amplitude $|u|$. Figure \ref{fig:sh_wave_surface} compares the MH-PINN predictions against the reference solutions computed by the Null-field Boundary Integral Equation Method (Null-field BIEM) for four different incident angles: grazing incidence ($\theta_{inc} = 0$), oblique incidences ($\theta_{inc} = \pi/6, \pi/3$), and vertical incidence ($\theta_{inc} = \pi/2$). The horizontal axis represents the projected $x$-coordinate of the ground and the canyon surface.

% --- 将峡谷散射的大图移到这里，即刚好在引出这段图的文字下方 ---
% 这样 LaTeX 会优先在这里排版图片，然后再排版后面对于物理现象的分析文字

The MH-PINN results exhibit excellent agreement with the reference numerical solutions across all incident angles. From a physical perspective, when the SH wave is incident obliquely or grazingly ($\theta_{inc} = 0, \pi/6, \pi/3$), intense wave fluctuations (peaks and valleys) are observed on the illuminated side ($x < 0$). This is caused by the multiple wave reflections and interference between the ground surface and the canyon topography. Conversely, in the region behind the canyon ($x > 0$), the displacement amplitude is significantly attenuated. This phenomenon is known as the shielding effect of a concave topography against oblique incident waves. As the incident angle $\theta_{inc}$ decreases, the shielding effect of the circular canyon becomes more pronounced. For the case of vertical incidence ($\theta_{inc} = \pi/2$), the shielding effect vanishes, and the displacement amplitude distribution is perfectly symmetric with respect to the origin.

These results compellingly demonstrate that the MH-PINN, empowered by the analytical compensation technique, can accurately and efficiently solve unbounded wave problems subjected to complex Neumann boundary conditions, overcoming the limitations of standard penalty-based PINNs.

\section{Discussion}
\label{sec:limitations}
The proposed MH-PINN successfully solves the challenges of truncation error and far-field radiation condition satisfaction in unbounded wave scattering problems by introducing analytical coordinate transformation and physical asymptotic fitting. Numerical experiments show that this framework exhibits significant accuracy and convergence robustness compared to traditional soft-constrained PINNs when dealing with high-frequency fluctuations and complex geometries. However, as an emerging computational framework, the application of MH-PINN to a broader range of physical scenarios reveals certain methodological limitations. This section provides a critical discussion of these challenges encountered in its current stage of development.

\subsection{Topological Adaptability of the Multi-Body Scattering Problem}

The existing MH-PINN framework relies on a global coordinate mapping strategy, which implicitly assumes that the scatterer must be simply connected. This single background coordinate transformation cannot be directly applied to multi-body scattering problems.

When multiple discrete scatterers exist in the physical domain, it is impossible to establish a unified coordinate system that simultaneously maps all object surfaces to a single boundary $\xi = -1$ of the computational domain. Forcing a single mapping center would lead to severe spatial distortion and topological overlap. To address this issue, future research will focus on combining MH-PINN with domain decomposition or overlapping mesh techniques. This involves establishing an independent locally mapped subnet for each scatterer and enabling information exchange between subnets in overlapping regions through interface conditions. This will allow MH-PINN to handle complex wave problems such as multi-body dynamic disturbances and metamaterial arrays.

\subsection{Coordinate Transformation Strategy and Oscillation Accumulation Problem}

There are other alternatives when dealing with unbounded domain problems, such as using a reciprocal mapping $s = 1/r$ to map the point at infinity to the origin. However, this strategy has fundamental mathematical flaws in wave problems, which is why this paper insists on using a specific algebraic mapping function instead of a simple reciprocal substitution.

Sound waves satisfy the Sommerfeld radiation condition at infinity, and its asymptotic solution is usually in the form of $u(r) \sim \frac{e^{\mathrm{i}kr}}{\sqrt{r}}$. If a reciprocal substitution $s = 1/r$ is used, the solution near the mapped domain becomes:
\begin{equation}
    u(s) \sim \sqrt{s} \cdot e^{\mathrm{i}\frac{k}{s}}, \quad s \rightarrow 0
    \label{eq:reciprocal_substitution}
\end{equation}

The phase factor $e^{\mathrm{i}\frac{k}{s}}$ reveals a serious numerical singularity: when $s \rightarrow 0$, the phase $\frac{k}{s}$ tends to infinity. This means that near the origin of the mapped domain, the wave function will experience an infinitely violent accumulation of oscillations. For neural networks, this means that the objective function has an infinitely high-frequency component at $s = 0$. Due to the inherent spectral bias of neural networks, it is extremely difficult for the network to capture this infinitely high-frequency oscillation feature.

Furthermore, the Sommerfeld radiation condition becomes:
\begin{equation}
    \lim_{s \to 0} \frac{1}{\sqrt{s}} \left( -s^2 \frac{\partial u}{\partial s} - \mathrm{i}ku \right) = 0
    \label{eq:sommerfeld_reciprocal}
\end{equation}
under the reciprocal mapping, introducing a singular term with respect to $s$, making the condition extremely difficult to enforce numerically. Therefore, simple spatial compaction is not suitable for wave equations, and a mapping strategy that can maintain waveform stationarity must be designed, as is the method adopted in this paper.

\subsection{Summary}

Despite the aforementioned limitations, MH-PINN's accurate satisfaction of both Dirichlet and Neumann boundary conditions, far-field radiation conditions, and geometric adaptability in single-body unbounded scattering problems demonstrate its immense potential as a robust intelligent wave solver. The current restrictions, primarily the topological extension to multi-body scattering, are not fundamental flaws of the framework. Rather, these challenges point to the need for further exploration, and incorporating domain decomposition techniques to address them will be the primary focus of our future work to drive the algorithmic evolution of MH-PINN.

\section{Conclusion}

To address the challenges in the numerical simulation of unbounded wave problems, such as large truncation boundary errors, high computational costs, and the inability of traditional deep learning methods to satisfy far-field radiation conditions, this paper proposes a Mapping-based Hard-constraint physics-informed neural network (MH-PINN). This method combines analytical coordinate transformation with a neural network architecture embedding physical priors to establish a high-precision computational framework that is mesh-free and requires no manual truncation boundaries. The main contributions and findings of this paper are summarized below:

\begin{itemize}
    \item The compaction coordinate mapping strategy effectively solves the sampling problem in unbounded domains. By analytically mapping the infinite physical domain to the finite parameter domain, MH-PINN not only avoids the spurious reflection problem caused by artificially truncated boundaries in traditional finite element or finite difference methods, but also achieves adaptive densification of collocation points in the near-field region, thus enabling more efficient capture of complex scattering features near objects or topographical irregularities.
    
    \item The physics-based hard-constraint architecture fundamentally improves the accuracy and stability of the solution. By explicitly constructing a solution structure in the network output layer that includes far-field asymptotic factors and inner-boundary distance weights, and introducing novel exact hard-constraint formulations for both Dirichlet and Neumann boundary conditions, this method ensures that the predicted solution strictly satisfies the given boundary conditions and Sommerfeld radiation conditions in its mathematical structure. This design completely eliminates the boundary penalty terms in the loss function, greatly reduces the optimization difficulty, and avoids the gradient competition problem in high-dimensional non-convex optimization.
    
    \item Numerical experiments demonstrate that MH-PINN exhibits exceptional robustness and accuracy across a diverse range of unbounded wave problems. It successfully models acoustic radiation and scattering over a wide frequency band (including strongly oscillatory cases up to $k=20$) and within inhomogeneous media. Furthermore, it exhibits strong geometric adaptability across complex two-dimensional and three-dimensional scatterers  and accurately resolves elastodynamic behaviors, such as SH wave scattering over a topographical canyon. Compared to standard penalty-based PINNs, MH-PINN significantly accelerates convergence and maintains stable, low relative errors, effectively overcoming the "spectral bias" common in deep learning solvers.
\end{itemize}

In summary, MH-PINN provides a universal, mesh-free, and highly efficient paradigm for solving partial differential equations in unbounded domains. Demonstrating excellent scalability across both two-dimensional and three-dimensional configurations, the framework establishes a solid foundation for complex wave field modeling. Future work will focus on integrating domain decomposition techniques to address multi-body scattering with complex topologies and extending the methodology to parameter inversion tasks, ultimately providing a powerful intelligent simulation tool for broad applications ranging from underwater acoustics and metamaterial design to seismological wave analysis.

\section*{Acknowledgments}
This work was supported by the National Natural Science Foundation of China (Grant Nos. 12402238, 12372196), the Open Fund Project of State Key Laboratory of Structural Analysis, Optimization and CAE Software for Industrial Equipment (Grant No. GZ25119), and the Fundamental Research Funds for the Central Universities (Grant No. B250201235).


%% 告诉 LaTeX 使用哈佛格式排版参考文献
\bibliographystyle{elsarticle-num} 

%% 告诉 LaTeX 去读取你文件夹里的 cas-refs.bib 文件
\bibliography{cas-refs}

@article{li2012analysis,
  title={Analysis of the scattering by an unbounded rough surface},
  author={Li, P and Shen, J},
  journal={Mathematical Methods in the Applied Sciences},
  volume={35},
  number={18},
  pages={2166--2184},
  year={2012}
}

@article{chandler2012numerical,
  title={Numerical-asymptotic boundary integral methods in high-frequency acoustic scattering},
  author={Chandler-Wilde, SN and Graham, IG and Langdon, S and Spence, EA},
  journal={Acta Numerica},
  volume={21},
  pages={89--305},
  year={2012}
}

@article{alves2024wave,
  title={Wave scattering problems in exterior domains with the method of fundamental solutions},
  author={Alves, CJS and Antunes, PRS},
  journal={Numerische Mathematik},
  volume={156},
  pages={375--394},
  year={2024}
}

@book{zienkiewicz1977finite,
  title={The Finite Element Method},
  author={Zienkiewicz, OE},
  year={1977},
  publisher={McGraw-Hill Inc},
  address={New York, NY}
}

@article{jiang2013numerical,
  title={Numerical solution of acoustic scattering by an adaptive DtN finite element method},
  author={Jiang, X and Li, P and Zheng, W},
  journal={Communications in Computational Physics},
  volume={13},
  number={5},
  pages={1277--1244},
  year={2013}
}

@article{li2015hybrid,
  title={Hybrid smoothed finite element method for acoustic problems},
  author={Li, E and He, ZC and Xu, X and others},
  journal={Computer Methods in Applied Mechanics and Engineering},
  volume={283},
  pages={664--688},
  year={2015}
}

@book{vanderEerden1996,
title = "Finite element method for acoustic radiation and scattering in an unbounded domain",
author = "van der Eerden, F.J.M.",
year = "1996",
language = "English",
series = "Aerospace Engineering Reports",
publisher = "National Aerospace Laboratory, NLR",
number = "NLR-TR 96603 L",
}

@book{rabczuk2019extended,
  title={Extended finite element and meshfree methods},
  author={Rabczuk, T and Song, JH and Zhuang, X and others},
  year={2019},
  publisher={Academic Press}
}

@article{zhang2023finite,
  title={Finite-Difference Frequency-Domain Scheme for Sound Scattering by a Vortex with Perfectly Matched Layers},
  author={Zhang, Y and Ling, Z and Du, H and others},
  journal={Mathematics},
  volume={11},
  number={18},
  pages={3959},
  year={2023}
}

@book{thomas2013numerical,
  title={Numerical partial differential equations: finite difference methods},
  author={Thomas, JW},
  year={2013},
  publisher={Springer Science \& Business Media}
}

@article{engquist1977absorbing,
  title={Absorbing boundary conditions for the numerical simulation of waves},
  author={Engquist, A and Majda, B},
  journal={Mathematics of Computation},
  volume={31},
  number={139},
  pages={629--651},
  year={1977}
}

@article{berenger1994perfectly,
  title={A perfectly matched layer for the absorption of electromagnetic waves},
  author={B{\'e}renger, J-P},
  journal={Journal of Computational Physics},
  volume={114},
  number={2},
  pages={185--200},
  year={1994}
}

@article{hao2022physics,
  title={Physics-informed machine learning: A survey on problems, methods and applications},
  author={Hao, Z and Liu, S and Zhang, Y and Ying, C and Feng, Y and Su, H and Zhu, J},
  journal={arXiv preprint arXiv:2211.08064},
  year={2022}
}

@article{wu1995direct,
  title={A direct boundary element method for acoustic radiation and scattering from mixed regular and thin bodies},
  author={Wu, TW},
  journal={The Journal of the Acoustical Society of America},
  volume={97},
  number={1},
  pages={84--91},
  year={1995}
}

@article{chandler2004high,
  title={A high-wavenumber boundary-element method for an acoustic scattering problem},
  author={Chandler-Wilde, SN and Langdon, S and Ritter, L},
  journal={Philosophical Transactions of the Royal Society of London. Series A: Mathematical, Physical and Engineering Sciences},
  volume={362},
  number={1816},
  pages={647--671},
  year={2004}
}

@article{liu2011recent,
  title={Recent advances and emerging applications of the boundary element method},
  author={Liu, YJ and Mukherjee, S and Nishimura, N and Schanz, M and Ye, W and Sutradhar, A and Pan, E and Dumont, NA and Frangi, A and Saez, A},
  journal={Applied Mechanics Reviews},
  volume={64},
  number={3},
  pages={030802},
  year={2011}
}

@article{xiu2007efficient,
  title={An efficient spectral method for acoustic scattering from rough surfaces},
  author={Xiu, D and Shen, J},
  journal={Communications in Computational Physics},
  volume={2},
  number={1},
  pages={54--72},
  year={2007}
}

@article{cho2019spectrally,
  title={Spectrally-accurate numerical method for acoustic scattering from doubly-periodic 3D multilayered media},
  author={Cho, MH},
  journal={Journal of Computational Physics},
  volume={393},
  pages={46--58},
  year={2019}
}

@article{shen2009some,
  title={Some recent advances on spectral methods for unbounded domains},
  author={Shen, J and Wang, L-L},
  journal={Communications in Computational Physics},
  volume={5},
  number={2--4},
  pages={195--241},
  year={2009}
}

@article{raissi2019physics,
  title={Physics-informed neural networks: A deep learning framework for solving forward and inverse problems involving nonlinear partial differential equations},
  author={Raissi, M and Perdikaris, P and Karniadakis, GE},
  journal={Journal of Computational Physics},
  volume={378},
  pages={686--707},
  year={2019}
}

@article{karniadakis2021physics,
  title={Physics-informed machine learning},
  author={Karniadakis, George Em and Kevrekidis, Ioannis G and Lu, Lu and Perdikaris, Paris and Wang, Sifan and Yang, Liu},
  journal={Nature Reviews Physics},
  volume={3},
  number={6},
  pages={422--440},
  year={2021}
}

@article{jeong2025advanced,
  title={An advanced physics-informed neural network-based framework for nonlinear and complex topology optimization},
  author={Jeong, H and Batuwatta-Gamage, C and Bai, J and Rathnayaka, C and Zhou, Y and Gu, Y},
  journal={Engineering Structures},
  volume={322},
  pages={119194},
  year={2025}
}

@article{zhongkai2024pinnacle,
  title={Pinnacle: A comprehensive benchmark of physics-informed neural networks for solving pdes},
  author={Zhongkai, H and Yao, J and Su, C and others},
  journal={Advances in Neural Information Processing Systems},
  volume={37},
  pages={76721--76774},
  year={2024}
}

@article{luo2025physics,
  title={Physics-informed neural networks for PDE problems: a comprehensive review},
  author={Luo, K and Zhao, J and Wang, Y and others},
  journal={Artificial Intelligence Review},
  volume={58},
  number={10},
  pages={323},
  year={2025}
}

@article{tang2022extrinsic,
  title={An extrinsic approach based on physics-informed neural networks for pdes on surfaces},
  author={Tang, Z and Fu, Z and Reutskiy, S},
  journal={Mathematics},
  volume={10},
  number={16},
  pages={2861},
  year={2022}
}

@article{fu2024physics,
  title={Physics-informed kernel function neural networks for solving partial differential equations},
  author={Fu, Z and Xu, W and Liu, S},
  journal={Neural Networks},
  volume={172},
  pages={106098},
  year={2024}
}

@article{xia2023spectrally,
  title={Spectrally adapted physics-informed neural networks for solving unbounded domain problems},
  author={Xia, M and B{\"o}ttcher, L and Chou, T},
  journal={Machine Learning: Science and Technology},
  volume={4},
  number={2},
  pages={025024},
  year={2023}
}

@article{ren2024seismicnet,
  title={SeismicNet: Physics-informed neural networks for seismic wave modeling in semi-infinite domain},
  author={Ren, P and Rao, C and Chen, S and others},
  journal={Computer Physics Communications},
  volume={295},
  pages={109010},
  year={2024}
}

@article{krishnapriyan2021characterizing,
  title={Characterizing possible failure modes in physics-informed neural networks},
  author={Krishnapriyan, Aditi and Gholami, Amir and Zhe, Shandian and Kirby, Robert and Mahoney, Michael W},
  journal={Advances in neural information processing systems},
  volume={34},
  pages={26548--26560},
  year={2021}
}

@inproceedings{rahaman2019spectral,
  title={On the spectral bias of neural networks},
  author={Rahaman, N and Baratin, A and Arpit, D and Draxler, F and Lin, M and Hamprecht, FA and Bengio, Y and Courville, A},
  booktitle={Proceedings of the 36th International Conference on Machine Learning (ICML)},
  pages={5301--5310},
  year={2019}
}

@article{wang2021eigenvector,
  title={On the eigenvector bias of Fourier feature networks: From regression to solving multi-scale PDEs with physics-informed neural networks},
  author={Wang, S and Wang, H and Perdikaris, P},
  journal={Computer Methods in Applied Mechanics and Engineering},
  volume={384},
  pages={113938},
  year={2021}
}

@article{wang2021understanding,
  title={Understanding and mitigating gradient flow pathologies in physics-informed neural networks},
  author={Wang, S and Teng, Y and Perdikaris, P},
  journal={SIAM Journal on Scientific Computing},
  volume={43},
  number={5},
  pages={A3055--A3081},
  year={2021}
}

@book{morse1986theoretical,
  title={Theoretical acoustics},
  author={Morse, Philip McCord and Ingard, K Uno},
  year={1986},
  publisher={Princeton university press}
}

@inproceedings{nair2024acoustic,
  title={Acoustic scattering simulations via physics-informed neural network},
  author={Nair, Siddharth and Walsh, Timothy F and Pickrell, Gregory and Semperlotti, Fabio},
  booktitle={Sensors and Smart Structures Technologies for Civil, Mechanical, and Aerospace Systems 2024},
  volume={12949},
  pages={138--145},
  year={2024},
  organization={SPIE}
}

@article{wang2022and,
  title={When and why PINNs fail to train: A neural tangent kernel perspective},
  author={Wang, Sifan and Yu, Xinling and Perdikaris, Paris},
  journal={Journal of Computational Physics},
  volume={449},
  pages={110768},
  year={2022},
  publisher={Elsevier}
}

\end{document}